\documentclass[11pt]{amsart}
\usepackage[foot]{amsaddr}
\usepackage{graphicx}
\usepackage[margin=1.25in]{geometry}
\usepackage[colorinlistoftodos]{todonotes}
\usepackage{amssymb}
\usepackage{amsmath}
\usepackage{amsthm}
\usepackage{amsfonts}
\usepackage{nicefrac}
\usepackage{mathtools}
\usepackage{tikz}
\usepackage{tikz-cd} 
\usepackage{verbatim}
\usepackage{mathabx}
\usepackage{listings}
\usepackage[colorlinks=true,citecolor=blue,linkcolor=blue,urlcolor=blue]{hyperref}
\usepackage{color}
\usepackage{bm}
\usepackage{enumitem}
\setlist[itemize]{leftmargin=*}
\usepackage[T1]{fontenc}
\usepackage{multirow}
\usepackage{booktabs}
\usepackage{adjustbox}
\usepackage{datetime}

\newdateformat{monthyeardate}{\monthname[\THEMONTH], \THEYEAR}
\definecolor{codegreen}{rgb}{0,0.6,0}
\definecolor{codegray}{rgb}{0.5,0.5,0.5}
\definecolor{codepurple}{rgb}{0.58,0,0.82}
\definecolor{backcolour}{rgb}{0.96,0.96,0.96}

\lstdefinestyle{mystyle}{
    backgroundcolor=\color{backcolour},   
    commentstyle=\color{codegreen},
    keywordstyle=\color{blue},
    numberstyle=\tiny\color{codegray},
    stringstyle=\color{codepurple},
    basicstyle=\ttfamily\footnotesize,
    breakatwhitespace=false,         
    breaklines=true,                 
    captionpos=b,                    
    keepspaces=true,                 
    numbers=left,                    
    numbersep=5pt,                  
    showspaces=false,                
    showstringspaces=false,
    showtabs=false,                  
    tabsize=2
}
\usetikzlibrary{arrows,calc,decorations,decorations.markings,fadings,positioning,patterns,shapes}
\tikzset{>=latex}
\tikzstyle{mypoint}=[inner sep=0pt,outer sep=0pt,minimum size=5pt,fill,circle]

\theoremstyle{definition}

\newtheorem*{theorem}{\textsc{\textbf{Theorem}}}          \newtheorem{theorem-N}[Thm]{\textsc{\textbf{Theorem}}}
              \newtheorem{lemma-N}[Thm]{\textsc{\textbf{Lemma}}}
      \newtheorem{corollary-N}[Thm]{\textsc{\textbf{Corollary}}}
 \newtheorem{proposition-N}[Thm]{\textsc{\textbf{Proposition}}}

    \newtheorem{definition-N}[Thm]{\textsc{\textbf{Definition}}}
\newtheorem*{remark}{\textsc{\textbf{Remark}}}            \newtheorem{remark-N}[Thm]{\textsc{\textbf{Remark}}}
          \newtheorem{example-N}[Thm]{\textsc{\textbf{Example}}}
      \newtheorem{observation-N}[Thm]{\textsc{\textbf{Observation}}}

\newtheorem{conjecture-N}[Thm]{\textsc{\textbf{Conjecture}}}

\newtheorem{algorithm-N}[Alg]{\textsc{\textbf{Algorithm}}}

\newenvironment{display}{\begin{center}\begin{tikzcd}}{\end{tikzcd}\end{center}}

\usepackage{etoolbox}
\newcommand{\define}[4]{\expandafter#1\csname#3#4\endcsname{#2{#4}}}

\forcsvlist{\define{\DeclareMathOperator}{}{}}{im,coker,rad,nil,Ann,Ass,codim,Spec,mSpec,diam,ord,Supp,supp,disc,Ob,vol,rank,Sym,Alt,Ind,,tr,spl,re, Jac, Pic, Div, gr, Nrd, Trd, Tr, Nm}
\forcsvlist{\define{\newcommand}{\mathrm}{}}{Hom,Mor,id,GL,SL,PSL,USp,PGL,SO,SU,U,Mat,Ext,Tor,Res,Cor,Inf,End,Irr,Aut,Gal,Lie,lcm,sign,triv,diag,Map,op,ev,act,alg,sep,unr,nr,ab,T,ad}

\forcsvlist{\define{\newcommand}{\mathsf}{}}{Set,Grp,Ab,CRing,Mod,Vect,Cat,Top,PreSh,Sh,Sch,Nat,Fun,Diff,Alg,Rep,Aff}

\forcsvlist{\define{\newcommand}{\mathbf}{}}{N,Z,Q,R,C,F,G,A,B,D}

\newcommand{\llrrparen}[1]{
  \left(\mkern -3mu\left( #1 \right)\mkern -3mu\right)}

\renewcommand{\O}{\mathcal{O}}
\renewcommand{\H}{\mathrm{H}}

\makeatletter
\newcommand*\bigcdot{\mathpalette\bigcdot@{.5}}
\newcommand*\bigcdot@[2]{\mathbin{\vcenter{\hbox{\scalebox{#2}{$\m@th#1\bullet$}}}}}
\makeatother

{\makeatletter
 \gdef\xxxmark{%
   \expandafter\ifx\csname @mpargs\endcsname\relax 
     \expandafter\ifx\csname @captype\endcsname\relax 
       \marginpar{xxx}
     \else
       xxx 
     \fi
   \else
     xxx 
   \fi}
 \gdef\xxx{\@ifnextchar[\xxx@lab\xxx@nolab}
 \long\gdef\xxx@lab[#1]#2{\textbf{[\xxxmark #2 ---{\sc #1}]}}
 \long\gdef\xxx@nolab#1{\textbf{[\xxxmark #1]}}
}

\usepackage{fancyhdr} 
\begin{document}

\title{Extensions of mod $p$ representations of division algebras over non-Archimedean local fields}

\author{Andrew Keisling}
\author{Dylan Pentland}

\begin{abstract} Let $F$ be a local field over $\Q_p$ or $\F_p\llrrparen{t}$, and let $D$ be a central simple division algebra over $F$ of degree $d$. In the $p$-adic case, we assume $p>de+1$ where $e$ is the ramification degree over $\Q_p$; otherwise, we need only assume $p$ and $d$ are coprime. For the subgroup $I_1=1+\varpi_D \O_D$ of $D^\times$ we determine the structure of $\H^1(I_1, \pi)$ as a representation of $D^\times/I_1$ for an arbitrary smooth irreducible representation $\pi$ of $D^\times$. We use this to compute the group $\Ext^1_{D^\times}(\pi,\pi')$ for arbitrary smooth irreducible representations $\pi$ and $\pi'$ of $D^\times$. In the $p$-adic case, via Poincar\'{e} duality we can compute the top cohomology groups and compute the highest degree extensions.
\end{abstract}

\maketitle

\tableofcontents

\setcounter{page}{1}

\section{Introduction} \label{sec:intro}

Let $F$ be a non-Archimedean local field whose residue field has characteristic $p$. In the local Langlands program, division algebras play an important role because of the Jacquet-Langlands correspondence. This correspondence states that for a division algebra $D$ of invariant $1/n$ over $F$ there is an injection
\begin{center}
\adjustbox{scale=1.4}{
\begin{tikzcd}
\left\{\substack{\text{irreducible smooth} \\ \text{representations} \\ \text{ of } D^\times \text{ over } \overline{\Q}_\ell} \right\}/\simeq  \ar[hook]{r} & \left\{ \substack{\text{irreducible smooth} \\ \text{representations}\\ \text{of }\GL_n(F) \text{ over } \overline{\Q}_\ell} \right\}/\simeq
\end{tikzcd}
}
\end{center}
where one can describe the image well. Here, we choose the prime $\ell$ so that $\ell\neq p$. The objects on the right side appear in the Langlands correspondence for $\GL_n$, which is then in natural bijection with certain Galois representations over $\overline{\Q}_\ell$.

When we instead take representations over a field of characteristic $p$, the problem becomes more difficult. There has been some progress towards proving a characteristic $p$ local Langlands correspondence, such as in \cite{breuil2010emerging}, but in general it is not clear how to formulate the Jacquet-Langlands correspondence given above. One attempt has been made in \cite{scholze2015p}, where for a division algebra of invariant $1/n$ over $F$, Scholze constructs a functor
\begin{center}
\adjustbox{scale=1.4}{
\begin{tikzcd}
\left\{ \substack{\text{smooth admissible} \\ \text{representations}\\ \text{of }\GL_n(F) \text{ over }\F_p} \right\} \ar{r}{\mathcal{F}} &  \left\{\substack{\text{smooth admissible} \\ \text{representations} \\ \text{ of } D^\times \text{ over } \F_p} \right\}
\end{tikzcd}
}
\end{center}
where $\mathcal{F}(\pi)$ also carries an action of $\Gal(\overline{F}/F)$. This gives some evidence for a mod $p$ analogue of the Jacquet-Langlands correspondence. Because of this, understanding the representation theory of $D^\times$ with $\overline{\F}_p$ coefficients can be expected to be of use in the mod $p$ local Langlands program. The irreducible representations of $D^\times$ over $\overline{\F}_p$ are easy to classify, however because $\Rep(D^\times)$ is not semisimple it is important to understand extensions between irreducible representations. These are what we will ultimately characterize.

Let $D$ be a central simple division algebra over $F$ of degree $d$, with a choice of uniformizer $\varpi_D$ and ring of integers $\O_D$. For $a|d$, set $D_a^\times = F^\times \O_D^\times \langle \varpi_D^a \rangle$. Then let $\pi = \Ind_{D_a^\times}^{D^\times}\chi$ and $\pi' = \Ind_{D_{a'}^\times}^{D^\times}\chi'$ be irreducible representations, where $\chi$ and $\chi'$ are characters of $D_a^\times$ and $D_{a'}^\times$. In \S\ref{sec:Ext1}, we prove the following main theorem with some mild constraints on the residue field characteristic of $F$.

\begin{theorem}[Theorem \ref{thm: main_thm}]
Let $\chi,\chi', \pi$ and $\pi'$ be as above. Let $S=\{\chi^s : s=\varpi_D^i, 0\le i < \gcd(a,a')\}$, so that the elements $s$ form a set of coset representatives for $D_{a'}^\times \setminus D^\times / D_a^\times$.

There are two types of direct summands in $\Ext^1_{D^\times}(\pi,\pi')$. If $\Res^{D_{a'}^\times}_{D_{\lcm(a,a')}^\times} \chi'$ is equal to some $\chi^s\in S$, we have a nonzero direct summand $A_{\chi^s}$ fitting into an exact sequence
\begin{display}
    0 \ar{r} & \overline{\F}_p \ar{r} & A_{\chi^s} \ar{r} & \H^1(1+\pi_F \O_F, \overline{\F}_p) \ar{r}& 0.
\end{display}
We also get a nonzero direct summand $A_{\chi^s} \simeq \overline{\F}_p$ for each $\chi^s\in S$ for which $\Res^{D_{a'}^\times}_{D_{\lcm(a,a')}^\times} \chi' \otimes (\chi^s)^\ast$ is extended trivially from a character $x \mapsto \left(\frac{x}{\sigma(x)}\right)^{p^i}$. Set $A_{\chi^s}=0$ otherwise. Then $\Ext^1_{D^\times}(\pi,\pi')\simeq \bigoplus_{\chi^s\in S} A_{\chi^s}$.
\end{theorem}
We note that the cohomology ring of the pro-$p$ group $1+\pi_F \O_F$ with the trivial action is known to be $\bigotimes_{i\in I} \overline{\F}_p[x_i]/x_i^2$ for some index set $I$ when $p$ is large enough so that there is no $p$-torsion. In this case $1+\pi_F\O_F$ is a direct product of $|I|$ copies of $\Z_p$ as a topological group. In the $p$-adic case we have $|I| = [F:\Q_p]$, and in the local function field case we have $I=\N$. Thus, $A_{\chi^s}$ of the first type described in the theorem has dimension $[F:\Q_p]+1$ as an $\overline{\F}_p$-vector space in the $p$-adic case and countable dimension in the local function field case.

This paper is organized as follows: in \S \ref{sec:prelims}, we go over some basic information about local division algebras and their representations. In \S\ref{sec:red to cohom}, we translate the problem of computing extensions to one of computing certain cohomology groups. Then, in \S \ref{sec:H1}, we compute the cohomology group $\H^1(1+\varpi_D \O_D, \pi)$ for any smooth irreducible representation $\pi$. Setting $I_1=1+\varpi_D \O_D$, we do this by computing the Frattini subgroup $[I_1,I_1]I_1^p$. We additionally show that almost all elements of this subgroup are of the form $[x,y]z^p$ for $x,y,z\in I_1$. More explicitly, when $d>4$ this is true for every element of the Frattini subgroup that does not also lie in the subgroup $I_3 = 1+\varpi_D^3\O_D$. In \S \ref{sec:Ext1}, we use this to compute $\Ext^1_{D^\times}(\pi,\pi')$ for arbitrary smooth irreducible representations $\pi$ and $\pi'$ of $D^\times$. In the $p$-adic case, we use Poincar\'{e} duality to compute higher extensions and also show how to get partial information about $\H^2(I_1,\overline{\F}_p)$. This is enough to compute all extension groups for a quaternion algebra over $\Q_p$.

\subsection{Notation} Throughout, we will fix some notation. We let $F$ be a $p$-adic field of ramification degree $e$ and residue field degree $f$ over $\Q_p$. We also allow $F$ to be an extension of $\F_p\llrrparen{t}$. We denote a choice of uniformizer by $\pi_F$, the ring of integers by $\O_F$, and the residue field by $k_F=\F_q$. The choice of discrete valuation is given by $\nu_F$, normalized so that $\nu_F(\pi_F) = 1$.

Over $F$, we consider a degree $d>1$ central simple division algebra $D$, so that $\dim_{F} D = d^2$. We use $\varpi_D$ for a choice of a uniformizer with respect to $\nu_D := \frac{1}{d}\nu_F \circ \Nrd$, where $\Nrd$ is the reduced norm on $D$. We moreover choose $\varpi_D$ so that $\varpi_D^d=\pi_F$. We denote the ring of integers by $\O_D$, and the residue field by $k_D=\F_{q^d}$. We fix an algebraic closure $\overline{\F}_p$ and embeddings $k_F\hookrightarrow k_D \hookrightarrow \overline{\F}_p$. The notation $[\cdot]$ is used for the Tiechmuller lift $[\cdot]: k_D \to D$. We let $\sigma$ be a generator of $\Gal(k_D/k_F)$ such that $\varpi_{D} [x] \varpi_D^{-1} = [\sigma(x)]$ for all $x\in k_D$, and set $I_1=1+\varpi_D \O_D$ to be the unique pro-$p$ Iwahori subgroup of $D^\times$. We assume that $p>de+1$ in the $p$-adic case, as this ensures there is no $p$-torsion in $I_1$. This is automatic in the local function field case; there we instead just assume $\gcd(p,d)=1$. For $a|d$ we set $D_a^\times := F^\times \O_D^\times \langle \varpi_D^a \rangle$.

Given a smooth representation $\pi$ of a locally profinite group $G$, we use $\Ext_G^n(\pi,-)$ to denote the $n$th derived functor of $\Hom_G(\pi,-)$, where $\Hom_G(\pi,\pi')$ is the space of $G$-equivariant linear maps from $\pi$ to another smooth $G$-representation $\pi'$. All representations will be over $\overline{\F}_p$. We use $\mathbf{1}$ to denote the trivial character of $G$, which is the space $\overline{\F}_p$ with the trivial action. If $\chi:G\to\overline{\F}_p^\times$ is a character of $G$, then we use $\chi^\ast$ to denote the dual character, so $\chi\otimes\chi^\ast = \mathbf{1}$. For a normal subgroup $H\le G$ and a character $\chi: H \to \overline{\F}_p^\times$, define $\chi^s(h)$ to be $\chi(s^{-1}hs)$ for $s\in G$ and $h\in H$ - we use this in the Mackey formula.

Throughout, the cohomology groups $\H^n(G,A)$ are continuous cohomology groups for a topological group $G$ and discrete $G$-module $A$. As a special case of this, $\Hom(G,A)$ is the group of continuous homomorphisms from $G$ to $A$, which we sometimes emphasize by writing it as $\Hom_{\mathrm{cont}}(G,A)$. For a compact open subgroup $K$ of a locally profinite group $G$, the Hecke algebra $\mathcal{H}_K$ is the algebra $\overline{\F}_p[K\setminus G/K]$ of $\overline{\F}_p$-valued bi-$K$-invariant continuous functions of compact support under convolution. For any group $G$ and $g,h\in G$, we use the convention $[g,h] = ghg^{-1}h^{-1}$. We use $[G,G]$ and $[G,G]G^p$ to denote the closures of the commutator and Frattini subgroups of an $\ell$-group group $G$. 

\section{Preliminaries} \label{sec:prelims}

\subsection{Structure of division algebras}

There is a decomposition
\[D^\times \simeq \varpi_D^\Z \ltimes (k_D^\times \ltimes I_1),\]
similar to the analogous decomposition of $F^\times$ except we must use semidirect products due to non-commutativity. The subgroup $I_1$ has a filtration
\[I_i = 1+\varpi_D^i \O_D,\]
where $I_i/I_{i+1} \simeq k_D$ as an additive group for $i\ge 1$. We may understand elements of $D$ explicitly via Tiechmuller lifts $[\cdot]: k_D \to D$, which allows a unique representation of every element as $x = \sum_{i\ge n} [x_i] \varpi_D^i$.

Up to isomorphism, the central simple division algebras $D$ over $F$ are classified by the result of the invariant map
\begin{display} 
\mathrm{Br}(F) \ar{r}{\mathrm{inv}_F} & \Q/\Z.
\end{display}
This is an isomorphism, and the isomorphism classes $[D]$ of degree $d$ central simple division algebras are obtained as preimages of $\frac{r}{d}$ where $\gcd(r,d) = 1$. The value $r$ reflects the choice of the generator $\sigma := (x\mapsto x^{q^r})$ of $\Gal(k_D/k_F)$, such that $\varpi_D[x]\varpi_D^{-1}=[\sigma(x)]$.

\subsection{Representations of \texorpdfstring{$D^\times$}{Dx}}

We work with smooth representations of $D^\times$ with coefficients in $\overline{\F}_p$, which makes the following well-known lemma crucial:

\begin{lemma-N}
Let $G$ be a pro-$p$ group. Then any nonzero smooth $\overline{\F}_p$-representation $V$ has $V^G \neq 0$.
\label{lem:pro-p_reps}
\end{lemma-N}

\begin{proof}
This is a classical result, which can be easily deduced from Proposition 26 in \cite{serre1977linear} which deals with finite $p$-groups.
\end{proof}

To begin, we would like to understand the characters of $D^\times$. There is an exact sequence
\begin{display}
1 \ar{r} & D^\times_{\Nrd=1} \ar{r} & D^\times \ar{r}{\Nrd} & F^\times \ar{r} & 1,
\end{display}
where exactness follows from the fact that $\Nrd$ is surjective. This is because the degree $d$ unramified extension $E/F$ is contained in $D$ and one can show $\Nrd|_{E}=\Nm_{E/F}$, so surjectivity of $\Nrd$ follows from surjectivity of the norm map from $\O_E^\times$ to $\O_F^\times$ for unramified extensions and the fact that $\Nrd(\varpi_D)$ is a unit multiple of $\pi_F$. It is known that $D^\times_{\Nrd=1} = [D^\times, D^\times]$, see for example \S 1.4.3 of \cite{platonov1993algebraic}. Then any character $D^\times \to \overline{\F}_p^\times$ arises as 
\begin{display}
\chi: D^\times \ar{r}{\Nrd} & F^\times \ar{r}{\kappa} & \overline{\F}_p^\times
\end{display}
because $\overline{\F}_p^\times$ is abelian, so any homomorphism into this group must factor through the abelianization. We can easily classify the characters $\kappa$ of $F^\times$ via $F^\times \simeq \pi_F^\Z \times k_F^\times \times (1+\pi_F \O_F)$, where the final component is pro-$p$ and hence we only need to compute characters of $\pi_F^\Z \times k_F^\times$ by the previous lemma. These are both cyclic groups, so the characters are then determined by where we send the generators of each component.

Now we turn to classifying characters of $D_a^\times$. Recall that this subgroup is defined to be $F^\times \O_D^\times \langle \varpi_D^a \rangle$ where $a|d$. As special cases, $D_d^\times = F^\times \O_D^\times$ and $D_1^\times=D^\times$. We have a decomposition
\[D_a^\times \simeq \varpi_D^{a\Z} \ltimes (k_D^\times \ltimes I_1),\]
and so we can denote elements by $(\varpi_D^{an},x,y)$ where $x\in k_D^\times$ and $y\in I_1$.

\begin{lemma-N}
The characters of $D_a^\times$ are given by
\[\chi_{a,\alpha,m}: (\varpi_D^{an},x,y) \mapsto \alpha^{n} \Nm_{k_D/\F_{q^a}}(x)^m.\]
Here, $\alpha\in \overline{\F}_p^\times$ and $m\in \Z$.
\end{lemma-N}

\begin{proof}
As $I_1$ is a pro-$p$ group, Lemma \ref{lem:pro-p_reps} implies that any character mapping to $\overline{\F}_p^\times$ must only depend on $\varpi_D^{an}$ and $x$. Because $\varpi_D^{a\Z}$ and $k_D^\times$ are cyclic and $k_D^\times$ embeds into $\overline{\F}_p^\times$, such a character must be of the form $\chi: (\varpi_D^{an},x,y) \mapsto \alpha^{n} x^m$. The only thing to work out is exactly which values of $m$ are allowed.

In order for $\chi$ to be well-defined on $D_a^\times$, it must be the case that
\[\chi([x])=\chi(\varpi_D^a [x] \varpi_D^{-a})=\chi([\sigma^a(x)]).\]
In particular, $x^m=\sigma^a(x^m)$ for all $x\in k_D^\times$. This gives a congruence condition on $m$. Namely, viewing $m \in \Z/|k_D^\times|\Z$, if $\sigma$ sends $x \mapsto x^{q^r}$ we have $m \equiv q^{ar}m \pmod{|k_D^\times|}$. This means $\frac{q^{d}-1}{q^{a}-1}$ divides $m$, and so the character restricted to $k_D^\times$ factors through the norm map $\Nm_{k_D/\F_{q^a}}$. From this we know that any character must take the form given in the lemma.
\end{proof}

\begin{corollary-N}
Suppose $a|a'$ are divisors of $d$. Then \[\Res^{D_a^\times}_{D_{a'}^\times} \chi_{a,\alpha,m} = \chi_{a',\alpha^{a'/a},m'}\]
where $m'=\frac{q^{a'}-1}{q^{a}-1}m$, which implies that $\Nm_{k_D/\F_{q^{a}}}(x)^m = \Nm_{k_D/\F_{q^{a'}}}(x)^{m'}$
for all $x\in k_D^\times$.
\label{cor:char_restriction}
\end{corollary-N}

In particular, we see that we can obtain all characters $\chi_{a',\alpha',0}$ from restrictions of characters $\chi_{1,\alpha,0}= \kappa \circ \Nrd$ of $D^\times$ as we may choose $\alpha$ to be an $a$th root of $\alpha'$ (recall $d$ is assumed to be coprime to $p$ in any case). We will now show that all irreducible representations can be obtained from inductions of characters of $D_a^\times$.

In the following theorem, a character $\chi$ of $k_D^\times$ is said to be of order $a$ if $a$ is the minimal integer such that $\chi(x)=\chi(\sigma^a(x))$ for all $x\in k_D^\times$. Such a character can be inflated along the quotient map $\O_D^\times \to k_D^\times$, and then extended trivially to $F^\times \O_D^\times$. By assigning the value of $\varpi_D^a$ to be $1$, it is then extended to $D_a^\times$.

\begin{theorem-N}[\cite{ly2013representations}, Prop. 1.3.1]
The smooth irreducible mod $p$ representations $V$ of $D^\times$ are given by 
\[V_{\chi,\kappa} := \Ind_{D_a^\times}^{D^\times} \left(\chi \otimes \kappa \right).\]
Here, $\kappa$ is a character of the type $\chi_{a,\alpha,0}$. The character $\chi$ is extended from an order $a$ character of $k_D^\times$.

We have $V_{\chi,\kappa} \simeq V_{\chi',\kappa'}$ if and only if $\kappa = \kappa'$ and $\chi = (\chi')^s$ for some $s\in D^\times$.
\label{thm:irrep_classification}
\end{theorem-N}

\begin{proof}
Let $V$ be an irreducible representation of $D^\times$ over $\overline{\F}_p$. Then, again by Lemma \ref{lem:pro-p_reps}, we know $V^{I_1}$ is a nonzero subrepresentation because $I_1$ is pro-$p$ and normal, so $V^{I_1}=V$. It follows that irreducible representations of $D^\times$ are in bijection with those of
\[D^\times/I_1 \simeq \varpi_D^\Z \ltimes k_D^\times.\]
The irreducible representations of this group over $\overline{\F}_p$ are well-known, as it is a semidirect product of $\varpi_D^\Z\simeq \Z$ and a finite abelian group with order prime to $p$. This is shown in \cite{serre1977linear} \S 8.2, and the irreducibles match up to those given in the theorem statement.
\end{proof}

\section{Reduction to cohomology}\label{sec:red to cohom}

We will begin by explaining how to reduce the computation of $\Ext^n_{D^\times}(\pi,\pi')$ to extensions between characters, after which we can reduce that problem to computation of certain cohomology groups. There are two main facts used in this reduction: first, that $\pi$ and $\pi'$ are induced from characters. Secondly, the inductions are from \textit{finite index} subgroups, so Frobenius reciprocity is a two-sided adjunction.

\begin{lemma-N}
For a locally profinite group $G$, let $H\le G$ be a closed subgroup such that $[G:H]<\infty$, and let $V$ and $W$ be smooth representations of the groups $H$ and $G$ respectively. Then we have 
\[\Ext^n_G(\Ind_H^G V, W) \simeq \Ext^n_H(V, \Res^G_H W).\]
This also holds in the other direction. 
\end{lemma-N}

\begin{proof}
This is a well-known fact from category theory: Frobenius reciprocity carries on to the derived functors of $\Hom$ because we have an adjunction both ways, since for finite index subgroups compact induction agrees with induction.
\end{proof}

Now we can use this to reduce our problem to computing extensions of characters. Consider irreducible representations $\pi = \Ind_{D_{a}^\times}^{D^\times}\chi$ and $\pi' = \Ind_{D_{a'}^\times}^{D^\times}\chi'$. Here, we have absorbed the $\kappa$ component in the notation of Theorem \ref{thm:irrep_classification} into $\chi$, so that $\chi$ and $\chi'$ are more general characters of $D_a^\times$ and $D_{a'}^\times$ rather than characters extended from $k_D^\times$. We would like to compute the dimension of $\Ext_{D^\times}^n(\pi,\pi')$. 


\begin{theorem-N}
Let $\chi,\chi',\pi$ and $\pi'$ be as given above. Then
\[\Ext^n_{D^\times}(\pi, \pi') \simeq \bigoplus_{\overline{s} \in D_{a'}^\times \setminus D^\times / D_a^\times} \Ext^n_{D_{\lcm(a,a')}^\times}(\mathbf{1}, (\Res^{D_{a'}^\times}_{D_{\lcm(a,a')}^\times} \chi' )\otimes (\chi^s)^\ast),\]
where $s$ is a coset representative of the double coset $\overline{s}$. The characters $\chi^s$ on $D_{\lcm(a,a')}^\times$ are defined as $\chi^s(x)=\chi(s^{-1}xs)$, and so are conjugated restrictions of $\chi$.
\label{thm:reduction_to_chars}
\end{theorem-N}

\begin{proof}
We first apply Frobenius reciprocity, the first time on the right induced representation. By the previous lemma, we may apply it for $\Ext^n$ as the subgroups $D_a^\times$ and $D_{a'}^\times$ have finite index in $D^\times$. We have
\[\Ext^n_{D^\times}(\pi,\pi') = \Ext^n_{D^\times} (\Ind_{D_a^\times}^{D^\times}\chi,\Ind_{D_{a'}^\times}^{D^\times}\chi') \simeq \Ext^n_{D_{a'}^\times}(\Res^{D^\times}_{D_{a'}^\times}\Ind_{D_{a}^\times}^{D^\times}\chi,\chi').\]
By the Mackey formula, we have
\[\Res^{D^\times}_{D_{a'}^\times}\Ind_{D_{a}^\times}^{D^\times}\chi \simeq \bigoplus_{\overline{s} \in D_{a'}^\times \setminus D^\times / D_{a}^\times} \Ind_{sD_a^\times s^{-1} \cap D_{a'}^\times}^{D_{a'}^\times} \chi^s,\]
where $s$ is a representative of the double coset $\overline{s}$. Because $D_a^\times$ is normal, the subgroup we induce from is simply $D_a^\times \cap D_{a'}^\times = D_{\lcm(a,a')}^\times$.

We now pull out the direct sum and apply Frobenius reciprocity on the other side. We then obtain
\[\bigoplus_{\overline{s} \in D_{a'}^\times \setminus D^\times / D_{a}^\times} \Ext^n_{D_{a'}^\times} (\Ind_{D_{\lcm(a,a')}^\times}^{D_{a'}^\times} \chi^s, \chi') \simeq \bigoplus_{\overline{s} \in D_{a'}^\times \setminus D^\times / D_{a}^\times} \Ext^n_{D_{\lcm(a,a')}^\times}(\chi^s, \Res^{D_{a'}^\times}_{D_{\lcm(a,a')}^\times} \chi').\]
The result then follows after tensoring with $(\chi^s)^\ast$ on each individual extension group.
\end{proof}

Thus, it suffices to be able to compute $\Ext^n_{D_a^\times}(\mathbf{1},\chi)$ for any character $\chi$ of $D_a^\times$.

We can now employ all the tools of group cohomology to solve our problem. It is tempting to try to use the Hochschild-Serre spectral sequence and identify $\Ext^n_{D_a^\times}(\mathbf{1},\chi)$ with $\H^n(D_a^\times, \chi)$, but this does not quite work. Since we must work with continuous cohomology and smooth representations, and $D_a^\times$ is not a profinite group, we cannot use the Hochschild-Serre spectral sequence in this case. However, we can recover a similar spectral sequence.

\begin{proposition-N}
There is a first quadrant spectral sequence
\[E_2^{i,j} = \H^i(D_a^\times/I_1,\H^j(I_1,\chi)) \implies \Ext^{i+j}_{D_a^\times}(\mathbf{1},\chi).\]
\label{prop: ext_spectral_seq}
\end{proposition-N}

\begin{proof}
From equation (33) in \S 9 of \cite{paskunas2007extensions}, letting $\mathcal{H}=\overline{\F}_p[D_a^\times/I_1]$ be the Hecke algebra of $I_1$, we have a spectral sequence
\[E_2^{i,j} = \Ext^i_{\mathcal{H}}(\mathbf{1},\H^j(I_1,\chi)) \implies \Ext^{i+j}_{D_a^\times}(\mathbf{1},\chi).\]
We may ignore the central character in our setting. The extension group is computed in the category of all $\mathcal{H}$-modules or equivalently $D_a^\times/I_1$-representations, which allows us to identify $\Ext^i_{D_a^\times/I_1, \mathrm{all}}(\mathbf{1},-) \simeq \H^i_{\mathrm{all}}(D_a^\times/I_1,-)$. Because the quotient $D_a^\times/I_1 \simeq \varpi_D^{a\Z} \ltimes k_D^\times$ is a discrete group, it follows that this is the same as the continuous cohomology. We then obtain a spectral sequence
\[E_2^{i,j} = \H^i(D_a^\times/I_1,\H^j(I_1,\chi)) \implies \Ext^{i+j}_{D_a^\times}(\mathbf{1},\chi)\]
as desired.
\end{proof}

As the spectral sequence of Proposition \ref{prop: ext_spectral_seq} is a first quadrant spectral sequence, we may obtain the following five-term exact sequence of low degree terms:
\begin{display}
    0 \ar{r} & \H^1(D_a^\times/I_1,\chi) \ar{r} & \Ext^1_{D_a^\times}(\mathbf{1},\chi) \ar{r}
             \ar[draw=none]{d}[name=X, anchor=center]{}
    & (\Hom(I_1,\overline{\F}_p) \otimes \chi)^{D_a^\times/I_1} \ar[rounded corners,
            to path={ -- ([xshift=2ex]\tikztostart.east)
                      |- (X.center) \tikztonodes
                      -| ([xshift=-2ex]\tikztotarget.west)
                      -- (\tikztotarget)}]{dll}[at end]{} \\      
    & \H^2(D_a^\times/I_1,\chi) \ar{r} & \Ext^2_{D_a^\times}(\mathbf{1},\chi)
\end{display}
Here, we have expanded the definition of $\H^1(I_1,\chi)^{D_a^\times/I_1}$. Thus, we obtain an exact sequence similar to the inflation-restriction exact sequence we would obtain via the Hochschild-Serre spectral sequence if $D_a^\times$ were a profinite group.

In Corollary \ref{cor: general_ext_sequence}, we use the same spectral sequence to obtain a description of $\Ext^n_{D_a^\times}(\mathbf{1},\chi)$ in terms of the cohomology of $I_1$, so this step in our method works for all degrees of extensions.

\section{Cohomology of \texorpdfstring{$I_1$}{I1}} \label{sec:H1}

Recall that we use $I_1$ to denote the subgroup $1+\varpi_D\O_D$ of $D^\times$. A key step in computing all extensions of irreducible representations of $D^\times$ is to understand the space $\H^1(I_1,\pi)$ as a representation of $D_a^\times/I_1$, where $\pi$ is a smooth irreducible representation of $D_a^\times$. Note that at the end of the previous section we had reduced to the case where $\pi$ was a character, but the computation is the same if $\pi$ is an irreducible of any dimension which is what we will assume here.

Since the action of $I_1$ is trivial on $\pi$, we have $\H^1(I_1,\pi) \simeq \Hom(I_1,\overline{\F}_p) \otimes \pi$, where the action of $g\in D_a^\times/I_1$ on $\Hom(I_1,\overline{\F}_p)$ sends a homomorphism $\varphi(x)$ to $\varphi(g^{-1}xg)$ and we view $\pi$ as a representation of $D_a^\times/I_1$. Note that we use the normality of $I_1$ here - in general, this cohomology group is only a module over the Hecke algebra $\mathcal{H}_{I_1}$. Thus, we seek to understand the homomorphism space $\Hom(I_1,\overline{\F}_p)$.

Because the additive group $\overline{\F}_p$ is abelian and every element is $p$-torsion, any homomorphism $\varphi:I_1\to \overline{\F}_p$ will factor through the following diagram:
\begin{display}
I_1 \ar{r}{\varphi} \ar{d} & \overline{\F}_p \\
\frac{I_1}{[I_1,I_1]I_1^p} \ar[dotted]{ru}{\varphi'} &
\end{display}
The homomorphism $\varphi'$ is unique given $\varphi$. The closed subgroup $[I_1,I_1]I_1^p\triangleleft I_1$ generated by commutators and $p$th powers is called the Frattini subgroup of $I_1$.

Hence, we can reduce our problem of computing $\Hom(I_1,\overline{\F}_p)$ to that of computing $\Hom(I_1/[I_1,I_1]I_1^p,\overline{\F}_p)$. We do this by computing the Frattini subgroup itself.

\subsection{The Frattini subgroup}

To compute $[I_1,I_1]I_1^p$, we will begin with computing $I_1^p$. In the $p$-adic case, this has a simple description.

\begin{proposition-N}
If $F$ is an extension of $\Q_p$, then $I_1^p = 1+\varpi_D^{de+1}\O_D$.
\label{prop:pth_powers}
\end{proposition-N}

\begin{proof}
First we show that $I_1^p\subset 1+\varpi_D^{de+1}\O_D$. Let $1+x\in I_1$ for $x\in \varpi_D\O_D$. Then we have $(1+x)^p = \sum_{0\le i\le p} \binom{p}{i}x^i = 1+px\sum_{1\le i\le p-1} \left(\binom{p}{i}/p\right)x^{i-1}+x^p$. Because $F$ is a $p$-adic field, we have $\varpi_D^{de} = \pi_F^e=pu$ for some unit $u\in \O_F^\times$. Thus, because $x\in \varpi_D\O_D$, $px$ is in $\varpi_D^{de+1}\O_D$. Moreover, because $p>de+1$, we have that $x^p$ is in $\varpi_D^{de+1}\O_D$, so $(1+x)^p$ is in $1+\varpi_D^{de+1}\O_D$ as desired.

Now we show that $I_1^p\supset 1+\varpi_D^{de+1}\O_D$. Let $1+y\in 1+\varpi_D^{de+1}\O_D$ for some $y\in \varpi_D^{de+1}\O_D$. As a formal power series, we know that
\[\left(\sum_{n=0}^\infty \binom{1/p}{n} Y^n\right)^p = 1+Y,\]
so to prove that $1+y\in I_1^p$ it suffices to show that $\sum_{n=0}^\infty \binom{1/p}{n} y^n$ converges in $I_1$.

Because the $n=0$ term in this series is $1$ and $\varpi_D\O_D$ is closed in $D$, it suffices to show that $\nu_D\left(\binom{1/p}{n} y^n\right)>0$ for all $n>0$ and that $\nu_D\left(\binom{1/p}{n} y^n\right)\to\infty$ as $n\to\infty$. We have the identity
\[\binom{1/p}{n} = \frac{(-1)^n}{n!} \prod_{0\le i \le n-1} \left(\frac{ip-1}{p}\right).\]
Thus, because $\nu_D = \nu_F = e\nu_p$ on $\Q_p$ (where $\nu_p$ is the $p$-adic valuation), we compute $\nu_D\left(\binom{1/p}{n} y^n\right) = \nu_D(y^n)-(e\nu_p(n!)+e\nu_p(p^n)) = n\nu_D(y)-e(\nu_p(n!)+n)$. We also have $\nu_D(\varpi_D) = \frac{1}{d}$, so because $y\in \varpi_D^{de+1}\O_D$, we have $\nu_D(y)\ge \nu_D(\varpi_D^{de+1}) = e\left(1+\frac{1}{de}\right)$. By Legendre's formula, $\nu_p(n!)$ is bounded above by $\frac{n}{p-1}$. Therefore,
\[\nu_D\left(\binom{1/p}{n} y^n\right)\ge ne\left(1+\frac{1}{de}\right)-ne\left(\frac{1}{p-1}+1\right),\]
and because $p-1>de$ this will be greater than $0$ for $n>0$ and approach $\infty$ as $n\to\infty$.
\end{proof}

When $F$ is an extension of $\F_p\llrrparen{t}$, we work in characteristic $p$, so we have $(1+x)^p = 1 + x^p$. This implies that
\[I_1^p = 1+(\varpi_D\O_D)^p.\]
We now turn to computing $[I_1,I_1]$. The following decomposition will be useful here, as well as later on.

\begin{lemma-N}
Suppose $p>de+1$ if $F/\Q_p$ or $\gcd(d,p)=1$ in the local function field case. We have a decomposition
\[I_1 = (1+\pi_F \O_F) \times I_{1,\Nrd=1},\]
where $I_{1,\Nrd=1}$ is the kernel of the reduced norm restricted to $I_1$.
\label{lem:subgrp_generators}
\end{lemma-N}

\begin{proof}
Because there is no $p$-torsion, $1+\pi_F \O_F$ is a direct product of copies of $\Z_p$ as a multiplicative group. It follows that every element has a unique $d$th root, which means that $\Nrd|_{1+\pi_F \O_F}: x \mapsto x^d$ is surjective. Any element $x\in I_1$ can then be divided by an element of $1+\pi_F \O_F$ to land in $I_{1,\Nrd=1}$, which shows $I_1 = (1+\pi_F \O_F)I_{1,\Nrd=1}$. These subgroups have trivial intersection, since the only $d$th root of $1$ in $1+\pi_F \O_F$ is $1$. Finally, both of these subgroups are normal because $1+\pi_F\O_F$ is central, so we have $I_1 = (1+\pi_F \O_F) \times I_{1,\Nrd=1}$ as desired.
\end{proof}

As $1+\pi_F \O_F$ is central in $I_1$, we then have $[I_1,I_1]=[I_{1,\Nrd=1},I_{1,\Nrd=1}]$. We compute this using \S 1.4 of \cite{platonov1993algebraic}. The following lemma is useful when considering the lowest nonzero coefficients in a commutator in $[I_{1,\Nrd=1},I_{1,\Nrd=1} \cap I_i]$. This will be used in the theorem that follows.

\begin{lemma-N}
Let $i\ge 0$, and for $y\in k_D$ let $\varphi_{i,y}\in \End_{k_F}(k_D)$ denote the map
\[\varphi_{i,y}: x\mapsto \sigma(x)y-x\sigma^i(y).\]
The image of $\varphi_{i,y}$ is the subspace
\[\ker (\Tr_{k_D/k_F}) \cdot \prod_{0\le j\le i} \sigma^j(y),\]
which has codimension one for $y\in k_D^\times$.
\label{lem:codim_1}
\end{lemma-N}

\begin{proof}
Since $\sigma$ fixes $k_F$, the $k_F$-linearity of this map follows immediately. When $y=0$, the result is clear, so suppose now that $y\in k_D^\times$.

We first compute the kernel of this map. We have $\varphi_{i,y}(x)=0$ if and only if $\sigma(x) y = x \sigma^i(y)$. For $x\neq 0$, this is equivalent to $\frac{x}{\sigma(x)} = \frac{y}{\sigma^i(y)}$, which always has a solution $x$ by the multiplicative version of Hilbert's theorem 90 as
\[\Nm_{k_D/k_F}(y/\sigma^i(y))=1\]
and $\sigma$ is a generator of the Galois group. Again since $\sigma$ is a generator, we have $k_D^{\langle\sigma\rangle} = k_F$, so the kernel is one-dimensional: supposing for nonzero $x'$ we have $\frac{x'}{\sigma(x')} = \frac{x}{\sigma(x)}$, then $\sigma(\frac{x'}{x})=\frac{x'}{x}$ so the ratio lies in $k_F^\times$. Thus, $\im \varphi_{i,y}$ is a codimension one subspace. Noting that
\[\frac{\varphi_{i,y}(x)}{\prod_{0\le j\le i} \sigma^j(y)} = \frac{\sigma(x)}{\prod_{0 < j \le i} \sigma^j(y)} - \frac{x}{\prod_{0\le j < i} \sigma^j(y)},\]
applying $\sigma$ to the second term on the right hand side yields the first term, so the right hand side is contained in $\ker(\Tr_{k_D/k_F})$ and
\[\im\varphi_{i,y}\subset\ker (\Tr_{k_D/k_F}) \cdot \prod_{0\le j\le i} \sigma^j(y).\]
The additive version of Hilbert's theorem 90 says $\ker \Tr_{k_D/k_F} = \im(\sigma(x)-x) = \im \varphi_{i,1}$, which was shown to also have codimension one. Hence, we have equality.
\end{proof}

\begin{theorem-N}
We have
\[[I_1,I_1]I_1^p=I_{1,\Nrd=1}I_1^p \cap I_2.\]
When $F/\Q_p$, we may describe $I_1^p$ as in Proposition \ref{prop:pth_powers}, and when $F/\F_p\llrrparen{t}$ we have $I_1^p=1+(\varpi_D \O_D)^p$. Here, we make the same assumptions on $p$ as Lemma \ref{lem:subgrp_generators}.
\label{thm:GL1_frattini_subgroup}
\end{theorem-N}

\begin{proof}
There is a filtration $I_{i,\Nrd=1} := I_{1,\Nrd=1} \cap I_i$. As previously discussed, by Lemma \ref{lem:subgrp_generators} we have
\[[I_1,I_1]=[I_{1,\Nrd=1},I_{1,\Nrd=1}].\]
By Theorem 1.9 in \cite{platonov1993algebraic}, for $d>2$ this commutator subgroup equals $I_{2,\Nrd=1}$. This can be extended to the case of $d=2$ here because of the restrictions on $p$. Because $I_1^p\subset I_2$, the claim of the theorem follows immediately.

The argument in \cite{platonov1993algebraic} is short, so we summarize it here. We show first that under the quotient map $q_i: I_i \to I_i/I_{i+1} \simeq k_D$ we have $q_{i+1}([I_{1,\Nrd=1}, I_{i,\Nrd=1}]) = q_{i+1}(I_{i+1,\Nrd=1})$. Indeed, a computation shows that the image $q_{i+1}([I_{1,\Nrd=1}, I_{i,\Nrd=1}])$ contains the elements $\varphi_{i,y}(x)$ where $y\in k_D$ and $x\in q_i(I_{i,\Nrd=1})$, where we have $q_i(I_{i,\Nrd=1}) = k_D$ if $d\nmid i$ and $q_i(I_{i,\Nrd=1}) = \ker(\Tr_{k_D/k_F})$ if $d|i$. Then Lemma \ref{lem:codim_1} can be used to show that these generate all of $q_{i+1}(I_{i+1,\Nrd=1})$ in either case: when $d|(i+1)$, we have $\im (\varphi_{i,y}) = \ker(\Tr_{k_D/k_F})$ for any $y$, so we indeed get all of $q_{i+1}(I_{i+1,\Nrd=1})$; when $d \nmid (i+1)$, we can show that we generate all of $q_{i+1}(I_{i+1,\Nrd=1})=k_D$ by first finding $y,y'\in k_D$ so that $\im (\varphi_{i,y}) \neq \im (\varphi_{i,y'})$. Choosing $y=1$ and $y'$ such that $\prod_{0\le j\le i} \sigma^j(y')$ does not lie in $k_F$, we can always do this. If, in addition, $d|i$, the values $x$ are restricted to $\ker( \Tr_{k_D/k_F})$ so we need to justify why this still suffices. The image of $\varphi_{i,y}|_{\ker (\Tr_{k_D/k_F})}$ still has codimension one since $k_F=\ker (\varphi_{i,y})$ has trivial intersection with $\ker( \Tr_{k_D/k_F})$ (because $p$ and $d$ are coprime) and so restricted to $\ker (\Tr_{k_D/k_F})$ the map $\varphi_{i,y}$ is an isomorphism. Thus, we have $q_{i+1}([I_{1,\Nrd=1}, I_{i,\Nrd=1}]) = q_{i+1}(I_{i+1,\Nrd=1})$.


It then follows that
\[[I_{1,\Nrd=1}, I_{i,\Nrd=1}] I_{i+2,\Nrd=1} = I_{i+1,\Nrd=1}.\]
Now $[I_{1,\Nrd=1},I_{1,\Nrd=1}]$ is a non-central normal subgroup of $D^\times_{\Nrd=1}$, and so by Theorem 3.3 in \cite{platonov1993algebraic} it is open and contains $I_{j,\Nrd=1}$ for some $j$. Supposing that the minimal such $j$ satisfies $j>2$, then we have
\[[I_{1,\Nrd=1},I_{1,\Nrd=1}] \supset [I_{1,\Nrd=1},I_{j-2,\Nrd=1}] I_{j,\Nrd=1} = I_{j-1,\Nrd=1},\]
which contradicts the fact that $j$ is minimal. Thus, we must have $j\le 2$. But we also have $[I_{1,\Nrd=1},I_{1,\Nrd=1}]\subset I_{2,\Nrd=1}$, so we have equality.
\end{proof}

\begin{remark}
In the case of $d=2$, as shown in \cite{riehm1970norm}, \S 5, the same result about the commutator subgroup of $I_{1,\Nrd=1}$ being $I_{2,\Nrd=1}$ will still hold so long as $D$ is not a dyadic division algebra. This case is ruled out because $\gcd(p,d)=1$ in any case, so this assumption is important to include.
\end{remark}

\subsection{A commutator construction}

While not needed for computing $\H^1(I_1,\pi)$, we can actually show that nearly all elements of $[I_1,I_1]I_1^p$ are products of a single commutator and a $p$th power, rather than products of many commutators and a $p$th power as done above. Namely, every element of $[I_1,I_1]I_1^p$ that does not lie in $I_3$ is a product of a commutator and $p$th power when $d>4$.

We first study $q_2([I_1,I_1])$ in detail in Proposition \ref{prop: comm_choice} to determine the exact number of ways a given value in $k_D$ can be produced. We use this in the corollary that follows to show there are enough ways to do this that we can inductively choose $x,y\in I_1$ so that $[x,y]$ produces the desired commutator by studying the coefficient of $\varpi_D^i$ in the expansion of $[x,y]$ using Tiechmuller representatives.

\begin{proposition-N}
Suppose $d>4$ and fix $\alpha \in k_D^\times$. There exist $1+[x]\varpi_D\neq 1$ and $1+[y]\varpi_D\neq 1$ such that \[q_2([1+[x]\varpi_D,1+[y]\varpi_D])=\alpha,\]
and $x/y$ does not lie in any proper subfield of $k_D$.
\label{prop: comm_choice}
\end{proposition-N}

\begin{proof}
If we show that there are more than $\frac{q^{d/2+1}-1}{q-1}$ ratios $x/y$ achieve any given value of $\alpha$, then the claim follows. This is because $\frac{q^{d/2+1}-1}{q-1}\ge \sum_{a\le d/2} q^a \ge |\bigcup_{a|d, a<d} \F_{q^a}|$, the number of elements in the union of all proper subfields of $k_D$. We first show that the number of ratios can be deduced from the number of pairs $(x,y)$.

To have $q_2(1+[x]\varpi_D,1+[y]\varpi_D])=\alpha$ requires $x \sigma(y)-\sigma(x)y=\alpha$. Viewed as a curve $C_\alpha$ over $k_D$, this has many symmetries. The $(q+1)$th roots of unity $\langle \zeta \rangle = \mu_{q+1}(k_D)$ act on $C_\alpha$. In particular, given $(x,y) \in C_\alpha(k_D)$ we have $(\zeta x, \zeta y)$ as another solution. If we fix the ratio $x/y$, then the system
\[x\sigma(y)-\sigma(x)y =\alpha, x/y = \beta\]
gives $x=y\beta$, so we solve $y\beta \sigma(y) - \sigma(y) \sigma(\beta)y = \alpha$. Note that $\beta-\sigma(\beta) \neq 0$ as $\alpha \neq 0$, so if solutions exist there are $|\mu_{q+1}(k_D)|$ of them. Thus, over nonzero values $\alpha$ the number of ratios $x/y$ from $(x,y)\in C_\alpha(k_D)$ is precisely $|C_\alpha(k_D)|/|\mu_{q+1}(k_D)|$, where $|C_\alpha(k_D)|$ denotes the number of points over $k_D$.

Thus, we now want to show $|C_\alpha(k_D)|$ exceeds $|\mu_{q+1}(k_D)|\frac{q^{d/2+1}-1}{q-1}$. If we are given $x,y\in k_D^\times$, then $(x,y)$ lies on some $C_\alpha(k_D)$ for $\alpha\in k_D$. Via Lemma \ref{lem:codim_1}, we may show exactly $|k_D^\times|\cdot |k_F^\times|$ of these yield $\alpha=0$. Then we conclude
\[\sum_{\alpha \in k_D^\times} |C_\alpha(k_D)| = (|k_D^\times|)^2-|k_D^\times|\cdot |k_F^\times|.\]
There are two relevant actions on the family of curves $C_\alpha$:
\begin{itemize}
    \item For $z\in k_D^\times$, if $(x,y)\in C_\alpha(k_D)$ then $(zx,zy)\in C_{z\sigma(z)\alpha}(k_D)$.
    \item For $A\in \GL_2(k_F)$, if $(x,y)\in C_\alpha(k_D)$ then $A \cdot (x,y)$ lies in $C_{\det A \cdot \alpha}(k_D)$.
\end{itemize}
Let $d$ be odd. The first point shows the cosets $k_D^\times/k_D^{\times(q+1)}$ take on common values as the image of $z\sigma(z)$ is identical to $z^{q+1}$. As $d$ is odd, the image of $z^{q+1}$ consists of all the squares. The second item shows that these two common values are the same. Thus, $|C_\alpha(k_D)|=|k_D|-|k_F|=q^d-q$ when $\alpha \in k_D^\times$ because we have produced bijections between $C_\alpha(k_D)$ over all $\alpha \in k_D^\times$.

Now let $d$ be even, and let $\alpha \in k_D$. We claim the curve $C_\alpha:x\sigma(y) - \sigma(x)y=\alpha$ has $k_D$ points in bijection with $X_{\alpha}:\sigma(x)x+\sigma(y)y=\zeta \alpha$. Here, we choose $\zeta \neq 0$ so that $\sigma(\zeta)=-\zeta$.

Explicitly, given $(x,y)\in X_{\alpha}$ we first choose $\omega$ such that $\sigma(\omega)\omega=-1$. Then $(x+\omega y, \omega x + y)$ is a solution to $x\sigma(y)+\sigma(x)y=(\sigma(\omega)+\omega)\zeta \alpha$. We have $(\sigma(\omega)+\omega)\in k_F^\times$. Dividing the $x$ coordinate by $\zeta$, we get a solution to $x\sigma(y)-\sigma(x)y=(\sigma(\omega)+\omega) \alpha$. Using the $\GL_2(k_F)$ action, the solutions are in bijection.

Finally, we observe that points in $X_{\alpha}(k_D)$ are in bijection with solutions to $x^{q+1}+y^{q+1}=\zeta \alpha$, simply because $x^{q+1}$ and $\sigma(x)x$ take on the same values and have equal numbers of preimages for each value. The explicit number of solutions to this general type of equation is computed in \cite{weil1949numbers}, and will satisfy the desired bound so long as $d > 4$. In particular, the values $|C_\alpha(k_D)|$ take on constant values for each coset of $k_D^\times/k_D^{\times(q+1)}$, with one coset for $\zeta \alpha = -1$ taking the value $q^{d}+1+(-q)^{d/2}\frac{q-q^3}{q+1}-(q+1)$ and the others a common value so that the average is $q^d-q$. The coset for $\zeta \alpha = -1$ counts points on $x^{q+1}+y^{q+1}+1=0$, which when projectivized (to add $q+1$ points) is the Fermat curve $x^{q+1}+y^{q+1}+z^{q+1}=0$. The zeta function of this curve is well-known over $\F_{q^2}$ - see \cite{shioda1979fermat}. That the remaining cosets take the same value can be seen from the explicit formula of \cite{weil1949numbers}. Thus, we have shown the desired bound on $|C_\alpha(k_D)|$ and the proposition follows.
\end{proof}

The following result is analogous to the result at the end of \S 1.4.3 of \cite{platonov1993algebraic}, which states that in $[D^\times, D^\times]$ every element is a product of at most two commutators (this was originally shown in \cite{wang1950commutator}).

\begin{corollary-N}
Suppose $d>4$. Any element $\alpha\in [I_1,I_1]I_1^p\setminus I_3$ can be constructed as a product of a commutator and a $p$th power.
\label{cor:comm_construct}
\end{corollary-N}

\begin{proof}
We wish to write $\alpha = [x,y] z^p$ for $x,y,z\in I_1$. Suppose we have $x = 1+\sum_{i\ge 1} [x_i]\varpi_D^i$, and similarly for $y$. We fix $x_1$ and $y_1$ to be as in Proposition \ref{prop: comm_choice}, so $q_2([x,y])=q_2(\alpha) \in k_D^\times$ and $x_1/y_1$ does not lie in a proper subfield of $k_D$. Let $i>1$ and assume $x_1, \ldots, x_{i-1}$ and $y_1, \ldots, y_{i-1}$ are fixed. A calculation shows that the $\varpi_D^{i+1}$ coefficient of $[x,y]$ is a Tiechmuller lift of
\[C + \varphi_{i,x_1}(y_i) + \varphi_{i,y_1}(x_i).\]
Here, $C\in k_D$ is a constant only depending on $x_1, \ldots, x_{i-1}$ and $y_1, \ldots, y_{i-1}$. To do this calculation, we recall that we have $[x]+[y]=[x+y]+O(p)$. If $F/\F_p\llrrparen{t}$, this error is zero, while in the $p$-adic case by Proposition \ref{prop:pth_powers} we can absorb it into the $I_1^p$ component $z^p$. We use this expression to construct commutators inductively.

By Lemma \ref{lem:codim_1}, we know $\{C + \varphi_{i,x_1}(y_i) + \varphi_{i,y_1}(x_i): x_i,y_i\in k_D\}\subset k_D$ is a sum of two codimension one subspaces, shifted by $C$. This is all of $k_D$ if and only if the two subspaces are distinct, which will depend on the choice of $x_1$ and $y_1$. By Lemma \ref{lem:codim_1}, $\im(\varphi_{i,x_1}) = \im(\varphi_{i,y_1})$ if and only if
\[\prod_{0\le j \le i} \sigma^j(x_1/y_1) \in k_F,\]
since $c \cdot \ker(\Tr_{k_D/k_F}) = \ker(\Tr_{k_D/k_F})$ if and only if $c\in k_F$ due to the non-degeneracy of the trace pairing. A calculation shows that this condition on $x_1/y_1$ holds if and only if $x_1/y_1$ is in the degree $\gcd(i+1,d)$ extension of $k_F$. Because $x_1/y_1$ does not lie in a proper subfield of $k_D$, we may obtain any Tiechmuller lift of $k_D$ as the coefficient of $\varpi_D^{i+1}$ of $[x,y]z^p$ when $d\nmid (i+1)$. When $d\mid (i+1)$, the values obtained for this coefficient are in the coset $C+\im(\varphi_{i,x_1})=C+\im(\varphi_{i,y_1})$, which has codimension one in $k_D$.

Using Lemma \ref{lem:subgrp_generators}, we can count the number of elements of $I_{1,\Nrd=1}$ modulo $I_i$ for any $i$ and use this to compute the number of elements of $[I_1,I_1]I_1^p\setminus I_3 = I_{1,\Nrd=1}I_1^p \cap (I_2\setminus I_3)$ modulo $I_i$. Doing this gives exactly the number of elements of $[I_1,I_1]I_1^p$ modulo $I_i$ that we can produce by multiplying a commutator and a $p$th power as above, thus proving the claim.
\end{proof}

Following through the same arguments by inspecting the values from Proposition \ref{prop: comm_choice} and comparing to $|\mu_{q+1}(k_D)| \frac{q^{d/2+1}-1}{q-1}$, we can get results for $d=2,3,4$ as well. For $d=2,3$ we can get the same final result, but for $d=4$ Proposition \ref{prop: comm_choice} does not hold since the exact number of points on $C_\alpha(k_D)$ need not exceed the bound $|\mu_{q+1}(k_D)| \frac{q^{d/2+1}-1}{q-1}$. However, it does hold for $\frac{q}{q+1}|k_D^\times|$ values in $k_D^\times$.

\subsection{Computing cohomology}

With our understanding of the Frattini subgroup, we can now compute the structure of $\H^1(I_1,\pi)$ as a representation of $D_a^\times/I_1$. For most of this section, we focus on computing the $D^\times/I_1$-representation structure where $\pi$ is a smooth irreducible representation of $D^\times$ because, as shown after the proof of Proposition \ref{prop:H1_module_struct}, we can obtain the $D_a^\times/I_1$-representation structure from this simply by restricting.

First we review why we have a $D^\times/I_1$ representation structure on $\H^1(I_1,\pi)$ for an irreducible smooth representation $\pi$. Being the pro-$p$ Iwahori subgroup of $D^\times$, there is the structure of a module over the Hecke algebra $\mathcal{H}_{I_1} \simeq \overline{\F}_p[I_1\setminus D^\times/I_1]$ on $\pi^{I_1}$. Since $I_1$ is normal, this has the structure of a representation of $D^\times/I_1$. The derived functors of this are the continuous cohomology groups $\H^i(I_1,\pi)$, hence there is a structure of a representation of $D^\times/I_1$ on $\H^1(I_1,\pi)$. As $\H^1(I_1,\pi) \simeq \H^1(I_1,\overline{\F}_p) \otimes \pi$, we will only need to compute this structure for the trivial representation. For $\H^1(I_1,\overline{\F}_p)$, the structure as a representation can be described very explicitly as the conjugation action of $D^\times/I_1$ on homomorphisms. Namely, $g\in D^\times/I_1$ sends a homomorphism $\varphi(x)$ to $\varphi(g^{-1}xg)$.

We now compute the homomorphisms in $\H^1(I_1,\overline{\F}_p)$ using Theorem \ref{thm:GL1_frattini_subgroup}.

\begin{theorem-N}
We have a decomposition
\[\H^1(I_1,\overline{\F}_p) \simeq \H^1(I_{1,\Nrd=1},\overline{\F}_p) \oplus \H^1(1+\pi_F \O_F, \overline{\F}_p),\]
where $\H^1(I_{1,\Nrd=1},\overline{\F}_p) \simeq \Hom(k_D,\overline{\F}_p)\simeq \overline{\F}_p^{df}$ and $\H^1(1+\pi_F \O_F, \overline{\F}_p)$ is isomorphic to $\overline{\F}_p^{ef}$ in the $p$-adic case and $\bigoplus_{i\in \N} \overline{\F}_p$ in the local function field case.
\label{thm:homomorphism_classification}
\end{theorem-N}

\begin{proof}
Due to Lemma \ref{lem:subgrp_generators}, the decomposition into these cohomology groups follows immediately. This is for the first cohomology group so we do not need more advanced methods, but it is worth noting this is a subcase of the K\"{u}nneth theorem for profinite groups.

By the proof of Theorem \ref{thm:GL1_frattini_subgroup}, we know that the Frattini subgroup of $I_{1,\Nrd=1}$ is $I_{1,\Nrd=1} \cap I_2 = I_{2,\Nrd=1}$. We then have
\[\H^1(I_{1,\Nrd=1}, \overline{\F}_p) \simeq \Hom(I_{1,\Nrd=1}/I_{2,\Nrd=1}, \overline{\F}_p).\]
The quotient $I_{1,\Nrd=1}/I_{2,\Nrd=1}$ is isomorphic to $k_D$: we know it is a finite dimensional $\F_p$-vector space since it is finite, abelian, and $p$-torsion. In the decomposition of Lemma \ref{lem:subgrp_generators}, upon taking the elements $1+[x]\varpi_D$ over $x\in k_D$, the $I_{1,\Nrd=1}$ components yield coset representatives of $I_{1,\Nrd=1}/I_{2,\Nrd=1} \subset k_D$ for any element of $k_D$. This is because $d\ge 2$, so division by an element of $1+\pi_F \O_F$ cannot affect the coefficient of $\varpi_D$. Thus, we have 
\[\H^1(I_{1,\Nrd=1},\overline{\F}_p) \simeq \Hom(k_D,\overline{\F}_p),\]
where this last space is isomorphic to $\overline{\F}_p^{df}$ because the additive group of $k_D$ is $\F_p^{df}$. Now we turn to $1+\pi_F \O_F$. In the $p$-adic case, the assumption $p>de+1$ ensures this group has no $p$-torsion, while in the local function field case this fact is a given. This implies that as a topological group we have $1+\pi_F \O_F$ is isomorphic to $\Z_p^{[F:\Q_p]}$ in the $p$-adic case, and $\Z_p^\N$ in the local function field case. As $\Hom_{\mathrm{cont}}(\Z_p,\overline{\F}_p) \simeq \overline{\F}_p$, in the $p$-adic case the claim follows. In the local function field case
\[\Hom_{\mathrm{cont}}(1+\pi_F \O_F, \overline{\F}_p) \simeq \Hom_{\mathrm{cont}}(\Z_p^\N, \overline{\F}_p) \simeq \bigoplus_{i\in \N} \Hom_{\mathrm{cont}}(\Z_p,\overline{\F}_p),\]
which is what we wanted.
\end{proof}

\begin{corollary-N}
For $0\le i < df$ and $\eta_i\coloneqq (x\mapsto x^{p^i}) \in \Aut_{\F_p}(k_D)$, define the homomorphisms $\phi^{\eta_i}:I_1\to \overline{\F}_p$ by the compositions
\begin{display}
I_1 \ar{r} & I_1/I_2 \ar{r}{\simeq} & k_D \ar{r}{\eta_i} & k_D \ar{r} & \overline{\F}_p
\end{display}
where the last map is the inclusion that we have fixed. These form a basis of the $\overline{\F}_p$-vector space $\H^1(I_{1,\Nrd=1},\overline{\F}_p) \subset \H^1(I_1,\overline{\F}_p)$.

For homomorphisms $\eta: 1+\pi_F \O_F \to \overline{\F}_p$, define homomorphisms $\psi^\eta:I_1\to\overline{\F}_p$ by
\begin{display}
I_1 \ar{r}{\Nrd} & 1+\pi_F\O_F \ar{r}{\eta} & \overline{\F}_p.
\end{display}
Taking a basis $\{\eta_j\}$ of $\Hom_{\mathrm{cont}}(1+\pi_F \O_F, \overline{\F}_p)$, the maps $\{\psi^{\eta_j}\}$ form a basis of $\H^1(1+\pi_F \O_F, \overline{\F}_p) \subset \H^1(I_1,\overline{\F}_p)$.
\label{cor:homs_basis}
\end{corollary-N}

\begin{proof}
We note that $\phi^{\eta_i}|_{1+\pi_F \O_F}=0$ and that $\psi^{\eta_j}|_{I_{1,\Nrd=1}}=0$, so the proposed basis elements lie in the correct components of the decomposition in Theorem \ref{thm:homomorphism_classification}, and we would like to show that they indeed form bases of these components.

We do this for the $\phi^{\eta_i}$ first. We already computed the dimension of $\H^1(I_{1,\Nrd=1},\overline{\F}_p)\subset \H^1(I_1,\overline{\F}_p)$, so it suffices to show that the $\phi^{\eta_i}$ are linearly independent. By surjectivity of the quotient map $I_1 \to I_1/I_2$, it suffices to show that the $\eta_i \in \Hom(k_D,\overline{\F}_p)$ are linearly independent. Viewing these as polynomials $x \mapsto x^{p^i} \in \overline{\F}_p[x]$, any linear combination that is equal to zero corresponds to a degree $< p^{df}$ polynomial having $p^{df}$ roots, which means the coefficients must all be zero. It follows that these homomorphisms are linearly independent.

Now consider the $\psi^{\eta_j}$. The reduced norm $\Nrd$ restricts to $x \mapsto x^d$ on $F$, so as in Lemma \ref{lem:subgrp_generators} we see $\Nrd:I_1\to 1+\pi_F\O_F$ is a surjection since $1+\pi_F \O_F \subset I_1$. The maps $\eta_j$ by definition give a basis of $\Hom_{\mathrm{cont}}(1+\pi_F \O_F, \overline{\F}_p)=\H^1(1+\pi_F \O_F, \overline{\F}_p)$, and so by surjectivity of the reduced norm the maps $\psi^{\eta_j}$ form a basis of $\H^1(1+\pi_F \O_F, \overline{\F}_p) \subset \H^1(I_1,\overline{\F}_p)$.
\end{proof}

We may now compute the structure of $\H^1(I_1,\overline{\F}_p)$ as a representation of $D^\times/I_1$.

\begin{proposition-N}
As a representation of $D^\times/I_1$ via the conjugation action, we have
\[\H^1(I_1,\overline{\F}_p) \simeq \H^1(1+\pi_F \O_F, \overline{\F}_p) \oplus \bigoplus_{i\in \Z/f\Z} \Ind_{D_d^\times/I_1}^{D^\times/I_1} \chi_{\eta_i}\]
where the action on $\H^1(1+\pi_F \O_F, \overline{\F}_p)$ is trivial and $\chi_{\eta_i}$ is extended trivially from the order $d$ character $k_D^\times\to\overline{\F}_p^\times$ given by $x \mapsto \left(\frac{\sigma(x)}{x}\right)^{p^i}$. Up to isomorphism, the choice of a coset representative of $i$ does not matter.
\label{prop:H1_module_struct}
\end{proposition-N}

\begin{proof}
By Corollary \ref{cor:homs_basis}, we see that every map in the component $\H^1(1+\pi_F \O_F, \overline{\F}_p) \subset \H^1(I_1,\overline{\F}_p)$ factors through $\Nrd$. As a result, the conjugation action of $D^\times/I_1$ is trivial on this component and it is a direct sum of copies of the trivial representation.

We now turn to the action on $\H^1(I_{1,\Nrd=1}, \overline{\F}_p)$. To do this, we will use the basis $\phi^{\eta_i}$ provided by Corollary \ref{cor:homs_basis}. As $D^\times/I_1 \simeq \varpi_D^\Z \ltimes k_D^\times$, it suffices to study the $k_D^\times$ and $\varpi_D^\Z$ actions. Supposing $\varpi_D [x] \varpi_D^{-1} = [\sigma(x)]$ and $\sigma$ sends $x \mapsto x^{q^r}$, the conjugation action of $\varpi_D$ is by
\[\varpi_D \cdot \phi^{\eta_i} = \phi^{\eta_{i-rf}}.\]
Similarly, $[x] \cdot \phi^{\eta_i} = \left(\frac{\sigma(x)}{x}\right)^{p^i} \phi^{\eta_i}$ for all $x\in k_D^\times$. Defining $V_i$ to be the span of the $\phi^{\eta_{i'}}$ where $i' \equiv i \pmod{f}$, each is a subrepresentation so we obtain
\[\H^1(I_{1,\Nrd=1}, \overline{\F}_p) = \bigoplus_{i\in \Z/f\Z} V_i.\]
Each $V_i$ is a dimension $d$ representation.  We have $\Res^{D^\times/I_1}_{D_d^\times/I_1} V_i \simeq \bigoplus_{i' \equiv i \pmod{f}} \chi_{\eta_{i'}}$ where we define $\chi_{\eta_{i'}}$ for $i'\in \Z/(df)\Z$ as sending $x\in k_D^\times$ to $\left(\frac{\sigma(x)}{x}\right)^{p^{i'}}$ and extend trivially to $D_d^\times/I_1$. It follows that \[\Hom_{D_d^\times/I_1}(\chi_{\eta_i}, \Res^{D^\times/I_1}_{D_d^\times/I_1}V_i)\neq 0,\]
which by Frobenius reciprocity means $\Hom_{D^\times/I_1}(\Ind_{D_d^\times/I_1}^{D^\times/I_1} \chi_{\eta_i}, V_i)\neq 0$. However, by the classification of irreducible representations, the induced representation is irreducible of dimension $d$. As $\dim V_i = d$, it follows that $V_i$ is isomorphic to $\Ind_{D_d^\times/I_1}^{D^\times/I_1} \chi_{\eta_i}$.

\end{proof}

For any smooth irreducible representation $\pi$ of $D^\times$, this immediately gives us 
\[\H^1(I_1,\pi) \simeq (\H^1(1+\pi_F \O_F, \overline{\F}_p) \otimes \pi) \oplus \bigoplus_{i\in \Z/f\Z} \Ind_{D_d^\times/I_1}^{D^\times/I_1}(\chi_{\eta_i} \otimes \Res_{D_d^\times/I_1}^{D^\times/I_1} \pi).\]
We have already determined how to compute the restriction of a character of $D^\times$ in general in Corollary \ref{cor:char_restriction}, so when $\pi$ is a character it is easy to write this down explicitly.

We can now derive the $D_a^\times$-representation structure for other values of $a|d$. For the trivial representation, this is
\[\Res^{D^\times/I_1}_{D_a^\times/I_1} \H^1(I_1,\overline{\F}_p) \simeq \H^1(1+\pi_F \O_F,\overline{\F}_p) \oplus \bigoplus_{i\in \Z/f\Z} \Res^{D^\times/I_1}_{D_a^\times/I_1} \Ind_{D_d^\times/I_1}^{D^\times/I_1} \chi_{\eta_i}.\]
The Mackey formula decomposes each term in the direct sum as \[\Res^{D^\times/I_1}_{D_a^\times/I_1} \Ind_{D_d^\times/I_1}^{D^\times/I_1} \chi_{\eta_i}\simeq \bigoplus_{\overline{s}\in D_a^\times \setminus D^\times/D_d^\times} \Ind_{D_d^\times/I_1}^{D_a^\times/I_1} \chi_{\eta_i}^s.\]
We can similarly compute the tensor product with a smooth irreducible representation $\pi$ of $D_a^\times$.

\begin{remark}
In \cite{koziol2018hecke}, for a split reductive $p$-adic group $G$ and for a pro-$p$ Iwahori subgroup $I_1$ of $G$, the group $\H^1(I_1,\overline{\F}_p)$ is computed as a module over the Hecke algebra $\mathcal{H}_{I_1}$. Our result is somewhat analogous to Theorem 6.4 in that paper in how $\H_1(I_1,\overline{\F}_p)$ splits into two components, but the fact that $D^\times$ is anisotropic modulo its center makes its structure far more interesting than that of $\GL_1(F)$. Notably, our decomposition as a module over $\mathcal{H}_{I_1}$ similarly does not involve supersingular modules, but computing the same cohomology group for $\GL_2(D)$ using Theorem \ref{thm:GL1_frattini_subgroup} does actually produce supersingular modules if we assume $d>1$, as would be expected from the main theorem of \cite{koziol2018hecke}.
\end{remark}

\section{Extensions of irreducibles} \label{sec:Ext1}

Having computed the $D_a^\times/I_1$-representation structure of $\H^1(I_1,\pi)$, we are now ready to use this to compute extensions.

\subsection{\texorpdfstring{$\Ext^1_{D^\times}$}{Ext1} for all irreducible representations}

Let $\chi$ be any smooth character of $D_a^\times$. To use the exact sequence of low degree terms at the end of \S\ref{sec:red to cohom} to compute $\Ext_{D_a^\times}^1(\mathbf{1},\chi)$, and ultimately to compute $\Ext_{D^\times}^1(\pi,\pi')$ for any $\pi$ and $\pi'$, we will also need to compute the structure of $\H^i(D_a^\times/I_1,\chi)$.

\begin{proposition-N}
We have $\H^i(D_a^\times/I_1, \chi)=0$ unless the action by $\chi:D_a^\times \to \overline{\F}_p^\times$ is trivial, in which case $\H^1(D_a^\times/I_1, \chi)\simeq \overline{\F}_p$ and $\H^i(D_a^\times/I_1, \chi)=0$ for $i>1$.
\label{prop:group_cohom_zero}
\end{proposition-N}

\begin{proof}
Because we have a semidirect product $D_a^\times/I_1 \simeq \varpi_D^{a\Z} \ltimes k_D^\times$, we have a normal subgroup $k_D^\times$ and a quotient $\varpi_D^{a\Z}$ making an exact sequence 
\begin{display}
0 \ar{r} & k_D^\times \ar{r} & D_a^\times/I_1 \ar{r} & \varpi_D^{a\Z} \ar{r} & 0.
\end{display}
The group in question is discrete, and hence without requiring cochains to be continuous the cohomology is the same as taking all cochains. We may therefore use the Hochschild-Serre spectral sequence in ordinary cohomology.

We first show $\H^{i}(k_D^\times, \chi)=0$ for any $\chi$ and $i>0$. Observe that the composition
\begin{display}
\H^i(k_D^\times,\chi) \ar{r}{\mathrm{res}} & \H^i(1,\chi) \ar{r}{\mathrm{cores}} & \H^i(k_D^\times, \chi)
\end{display}
is multiplication by $|k_D^\times|$, and hence an isomorphism as we work in characteristic $p$. But the middle group is $0$, so this is also the zero map. It follows the cohomology group $\H^i(k_D^\times, \chi)$ is trivial. One can also see this by observing that $\widehat{\H}^0(k_D^\times,\chi)$ and $\widehat{\H}_0(k_D^\times,\chi)$ are both $0$ since $|k_D^\times|$ is prime to $p$, and then using Tate periodicity.

Because $\H^{i}(k_D^\times, \chi)=0$ is trivial for $i>0$, $\H^i(\varpi_D^{a\Z}, \chi^{k_D^\times}) \simeq \H^i(\varpi_D^{a\Z} \ltimes k_D^\times, \chi)$ via the inflation map in the higher inflation-restriction exact sequences arising from the Hochschild-Serre spectral sequence. The cohomological dimension of $\varpi_D^{a\Z}\simeq\Z$ is $1$, so for $i>1$ we have now shown the claim. For $i=1$, $\H^1(\varpi_D^{a\Z},\chi^{k_D^\times})=0$ unless $\chi^{k_D^\times}$ is trivial, in which case we get $\overline{\F}_p$. We only get $\chi^{k_D^\times}$ as the trivial $\varpi_D^{a\Z}$-character when $\chi$ itself is trivial.
\end{proof}

With this, we can now determine $\Ext^1_{D_a^\times}(\mathbf{1},\chi)$ via the spectral sequence Proposition \ref{prop: ext_spectral_seq}. We also obtain the following corollary:

\begin{corollary-N}
There is an exact sequence
\begin{display}
0 \ar{r} & \H^1(D_a^\times/I_1, \H^{i-1}(I_1,\chi)) \ar{r} & \Ext_{D_a^\times}^{i}(\mathbf{1},\chi) \ar{r} &  \H^{i}(I_1,\chi)^{D_a^\times/I_1} \ar{r} &  0.
\end{display}
\label{cor: general_ext_sequence}
\end{corollary-N}

\begin{proof}
Proposition \ref{prop:group_cohom_zero} shows that $\H^i(D_a^\times/I_1,V)=0$ for $i\ge 2$ for any representation $V$ of $D_a^\times/I_1$. Indeed, the previous argument can be adapted directly for $i\ge 2$ because $\Z$ still has cohomological dimension $1$ and $\Rep(k_D^\times)$ is semisimple as $|k_D^\times|$ is prime to $p$, and so we can decompose $V = \bigoplus_j \chi_j$ as a representation of $k_D^\times$ to get $\H^i(k_D^\times,V)=\bigoplus_j \H^i(k_D^\times,\chi_j) = 0$. We have a spectral sequence 
\[E_2^{i,j} = \H^i(D_a^\times/I_1,\H^j(I_1,\chi)) \implies \Ext^{i+j}_{D_a^\times}(\mathbf{1},\chi)\]
via Proposition \ref{prop: ext_spectral_seq} which we now know consists of two columns on the $E_2$ page, from which we get the above exact sequence as the $E_2$ page equals the $E_\infty$ page.
\end{proof}

We now compute $\Ext^1_{D_a^\times}(\mathbf{1},\chi)$.

\begin{theorem-N}
Let $D$ be a degree $d$ division algebra over $F$. Let $\chi$ be a character of $D_a^\times$ where $a|d$. There are two cases where the extensions $\Ext^1_{D_a^\times}(\mathbf{1},\chi)$ can be nontrivial:
\begin{itemize}
    \item When $\chi$ is trivial, there is an exact sequence
    \begin{display}
    0 \ar{r} & \overline{\F}_p \ar{r} & \Ext^1_{D_a^\times}(\mathbf{1},\chi) \ar{r} & \H^1(1+\pi_F \O_F,\overline{\F}_p) \ar{r}& 0
    \end{display}
    where $\H^1(1+\pi_F \O_F,\overline{\F}_p) \simeq \Hom(I_1, \overline{\F}_p)^{D_a^\times/I_1}$ is as in Theorem \ref{thm:homomorphism_classification}.
    
    \item When $a=d$, and $\chi$ is extended trivially from a character $x \mapsto \left(\frac{x}{\sigma(x)}\right)^{p^i}$ of $k_D^\times$ for some $i$, we have $ \Ext^1_{D_a^\times}(\mathbf{1},\chi) \simeq \overline{\F}_p$.
\end{itemize}
Otherwise, $\Ext^1_{D_a^\times}(\mathbf{1},\chi)=0$.
\label{thm:char_exts}
\end{theorem-N}

\begin{proof}
Due to Proposition \ref{prop:group_cohom_zero}, in the exact sequence of low degree terms arising from Proposition \ref{prop: ext_spectral_seq}, we obtain
\begin{display}
0 \ar{r} & \H^1(D_a^\times/I_1,\chi) \ar{r} & \Ext^1_{D_a^\times}(\mathbf{1},\chi) \ar{r} & (\Hom(I_1, \overline{\F}_p) \otimes \chi)^{D_a^\times/I_1} \ar{r} & 0.
\end{display}
When $\chi = \mathbf{1}$, by Proposition \ref{prop:H1_module_struct} we see that $\Hom(I_1, \overline{\F}_p)^{D_a^\times/I_1} \simeq \H^1(1+\pi_F \O_F, \overline{\F}_p)$ since this is the trivial component of the representation. Additionally, Proposition \ref{prop:group_cohom_zero} tells us that $\H^1(D_a^\times/I_1,\overline{\F}_p)\simeq \overline{\F}_p$. We then recover the first case of the theorem statement.

When $\chi$ is nontrivial, $\H^1(D_a^\times/I_1,\chi) = 0$ so we have $\Ext^1_{D_a^\times}(\mathbf{1},\chi) \simeq (\Hom(I_1, \overline{\F}_p) \otimes \chi)^{D_a^\times/I_1}$ via the restriction map. We know this as a $D_a^\times/I_1$-representation, so we just need to compute the trivial component.

Recall that as a $D_a^\times/I_1$-representation, we have already shown
\[\H^1(I_1,\chi) \simeq (\H^1(1+\pi_F \O_F, \overline{\F}_p)\otimes \chi) \oplus \bigoplus_{i\in \Z/f\Z} (\Res^{D^\times}_{D_a^\times} \Ind_{D_d^\times}^{D^\times} \chi_{\eta_i}) \otimes \chi.\]
As $\H^1(1+\pi_F \O_F, \overline{\F}_p)$ has a trivial action and $\chi$ is nontrivial, we know the $(\H^1(1+\pi_F \O_F, \overline{\F}_p) \otimes \chi)$ component is nontrivial. By the Mackey formula, the remaining component before tensoring with $\chi$ is
\[\bigoplus_{i\in \Z/f\Z} \Res^{D^\times}_{D_a^\times} \Ind_{D_d^\times}^{D^\times} \chi_{\eta_i} \simeq \bigoplus_{i\in \Z/f\Z} ~ \bigoplus_{\overline{s} \in D_a^\times \setminus D^\times / D_d^\times} \Ind_{D_d^\times}^{D_a^\times} \chi^s_{\eta_i}.\]

By the Mackey irreducibility criterion, $\Ind_{D_d^\times}^{D_a^\times} \chi^s_{\eta_i}$ is irreducible as it breaks down as a direct sum of distinct characters upon restriction to $D_d^\times$.

When we tensor with $\chi$, we can further pull $\chi$ into the induction via the push-pull formula to obtain a direct sum of inductions of the form
\[\Ind_{D_d^\times}^{D_a^\times} (\chi_{\eta_i}^s\otimes\Res_{D_d^\times}^{D_a^\times}\chi) \eqqcolon \Ind_{D_d^\times}^{D_a^\times} \chi'.\]
This remains irreducible since we have tensored with a chraracter. Now the trivial representation is irreducible as well, so if $\Ind_{D_d^\times}^{D_a^\times} \chi'$ contains a copy of the trivial representation then it must itself be trivial. It follows that we must have $a=d$ in this case, and also that $\chi_{\eta_i}^s\otimes\Res_{D_d^\times}^{D_a^\times}\chi \simeq \mathbf{1}$. As $a=d$, this just says that $\chi$ must be dual to some $\chi_{\eta_i}^s$, which can only be the case for at most one $i$ and $\overline{s}$. This is precisely the criterion in the theorem.
\end{proof}

This gives a way to compute any particular extension group.

\begin{theorem-N}
Let $\chi,\chi', \pi$ and $\pi'$ be as in Theorem \ref{thm:reduction_to_chars}. Let $S=\{\chi^s : s=\varpi_D^i, 0\le i < \gcd(a,a')\}$, so that the elements $s$ form a set of coset representatives for $D_{a'}^\times \setminus D^\times / D_{a}^\times$.

There are two types of direct summands in $\Ext^1_{D^\times}(\pi,\pi')$. If $\Res^{D_{a'}^\times}_{D_{\lcm(a,a')}^\times} \chi'$ is equal to some $\chi^s\in S$, we have a nonzero direct summand $A_{\chi^s}$ fitting into an exact sequence
\begin{display}
    0 \ar{r} & \overline{\F}_p \ar{r} & A_{\chi^s} \ar{r} & \H^1(1+\pi_F \O_F, \overline{\F}_p) \ar{r}& 0.
\end{display}
We also get a nonzero direct summand $A_{\chi^s} \simeq \overline{\F}_p$ for each $\chi^s\in S$ for which $\Res^{D_{a'}^\times}_{D_{\lcm(a,a')}^\times} \chi' \otimes (\chi^s)^\ast$ is extended trivially from a character $x \mapsto \left(\frac{x}{\sigma(x)}\right)^{p^i}$.

Set $A_{\chi^s}=0$ otherwise. Then $\Ext^1_{D^\times}(\pi,\pi')\simeq \bigoplus_{\chi^s\in S} A_{\chi^s}$.
\label{thm: main_thm}
\end{theorem-N}

\begin{proof}
We get this result simply by applying the result of Theorem \ref{thm:char_exts} in Theorem \ref{thm:reduction_to_chars}.
\end{proof}

If $\lcm(a,a')<d$, only the first type of direct summand occurs a single time; it is only in the case that $\lcm(a,a')=d$ that the situation can be more complicated.

We illustrate the example of $d=2$ below.

\begin{example-N}\label{ex:ext_1_quat_alg}
Let $D$ be a quaternion algebra over a $p$-adic field $F$. We will consider the case where $\pi=\Ind_{D_2^\times}^{D^\times} (\chi \otimes \kappa)$ and $\pi'=\Ind_{D_2^\times}^{D^\times} (\chi' \otimes \kappa')$ (in the notation of Theorem \ref{thm:irrep_classification}), since this is the case where the most interesting types of extensions occur.

Explicitly, these characters of $D_2^\times$ are of the form $\chi_{2,\alpha,m}$, where $q+1 \nmid m$ so that the characters have order $2$. There are four different possibilities for $\dim \Ext^1_{D^\times}(\pi,\pi'): 0, 1, ef+1,$ and $ef+2$.

We obtain dimension $ef+1$ when $\kappa=\kappa'$ and $\chi \otimes (\chi')^{\ast}$ is trivial, or $\kappa=\kappa'$ and $\chi \otimes ({\chi'}^{\varpi_D})^\ast$ is trivial. The character ${\chi'}^{\varpi_D}$, when we set $\chi'=\chi_{2,1,m'}$, is $\chi_{2,1,qm'}$. These conditions cannot both occur at the same time, which is why dimension $2ef+2$ is not possible.

Dimension $ef+2$ can occur when one of the above conditions holds, for example $\kappa=\kappa'$ and $\chi \otimes (\chi')^{\ast}$ is trivial. We then additionally require $\chi \otimes ({\chi'}^{\varpi_D})^\ast$ to be some character $x \mapsto (x^{1-q})^{p^i}$ when restricted to $k_D^\times$. This situation can occur, as well as the opposite situation where $\kappa=\kappa'$ and $\chi \otimes ({\chi'}^{\varpi_D})^\ast$ is trivial and $\chi \otimes ({\chi'})^\ast$ is some character $x \mapsto (x^{1-q})^{p^i}$ on $k_D^\times$.

We can obtain dimension $1$ when $\chi \otimes (\chi')^{\ast}$ is $x \mapsto (x^{1-q})^{p^i}$ on $k_D^\times$, but $\chi \otimes ({\chi'}^{\varpi_D})^\ast$ is nontrivial and vice versa. Finally, we get dimension $0$ in all other cases.

\end{example-N}

\subsection{Higher degree extensions in the \texorpdfstring{$p$}{p}-adic case}

In \cite{lazard1965groupes}, many fundamental results about $p$-adic analytic pro-$p$ groups are shown which will be helpful in understanding the structure of the higher cohomology groups in the case where $F$ is a $p$-adic field. Here, by a $p$-adic analytic group we mean a topological group with the structure of a $p$-adic analytic manifold over $\Q_p$ such that the group multiplication and inversion operations are analytic.

\begin{theorem}[\cite{lazard1965groupes} \S 2.5.8]
Let $G$ be a $p$-adic analytic pro-$p$ group and with no $p$-torsion, and let $r=\dim_{\Q_p} G$. Then $G$ is a Poincar\'{e} duality group of dimension $r$ over $\F_p$. That is,
\begin{itemize}
    \item $\H^n(G,\F_p)$ is finite dimensional for all $n\in \N$.
    \item $\dim_{\F_p} \H^r(G,\F_p)=1$.
    \item The cup product
    \[\H^n(G,\F_p) \times \H^{r-n}(G,\F_p) \to \H^r(G,\F_p)\]
    is a non-degenerate bilinear form.
\end{itemize}
Here, as always, all cohomology taken is continuous.
\end{theorem}

We will be interested in the case of $I_1=G$, and we will also want to understand how the cup product interacts with the $D_a^\times/I_1$-representation structure on the cohomology group. We note that $I_1$ is a $p$-adic analytic group as it is an open subgroup of $D^\times$, and since $p>de+1$ it has no $p$-torsion. Thus, $I_1$ is a Poincar\'{e} duality group.

The most difficult part of the following result is that the $D_a^\times/I_1$-representation structure on the top cohomology is trivial. This is the main result of \cite{koziol2020functorial} for general connected reductive groups. Using the method of Theorem 7.2 in \cite{koziol2018hecke}, we can alternatively show this by finding the uniform pro-$p$ subgroup $I_{de+1}$ and computing explicitly the action on its first cohomology, which determines the action on its top cohomology as the cohomology ring is an exterior algebra. We can then show that this action is trivial, and that it must agree with the action on the top cohomology of $I_1$ using the corestriction map.

\begin{proposition-N}[\cite{koziol2020functorial}]
\label{prop:poincare_duality}
Let $r=d^2ef$. As $D_a^\times/I_1$-representations, we have that $\H^n(I_1,\overline{\F}_p)^\ast \simeq \H^{r-n}(I_1,\overline{\F}_p)$. Here, $V^\ast$ denotes the dual representation.
\end{proposition-N}

\begin{proof}
As representations of $D^\times/I_1$ are equivalent to modules over $\mathcal{H}_{I_1}$ and $\dim_{\Q_p} I_1=d^2ef$, the structure on $\H^r(I_1,\F_p)$ as a $D^\times/I_1$-representation, and hence as a $D_a^\times/I_1$-representation, is trivial by \cite{koziol2020functorial}. The cup product
\[\H^n(I_1,\F_p) \times \H^{r-n}(I_1,\F_p) \to \H^r(I_1,\F_p) \simeq \F_p\]
is a non-degenerate bilinear form, which behaves well with respect to the $D_a^\times/I_1$-action. Let $d$ be an element of $D_a^\times/I_1$. Then we have for $\alpha \in \H^n(I_1,\F_p)$ and $\beta \in \H^{r-n}(I_1,\F_p)$ that
\[d \cdot (\alpha \smile \beta) = (d\cdot \alpha) \smile (d\cdot \beta).\]
It follows from the triviality of the top cohomology that Poincar\'{e} duality yields an isomorphism $\H^n(I_1,\F_p)^\ast \simeq \H^{r-n}(I_1,\F_p)$. Tensoring with $\overline{\F}_p$ gives the desired isomorphism.
\end{proof}

We will now again use the spectral sequence
\[E_2^{i,j} = \H^i(D_a^\times/I_1,\H^j(I_1,\chi)) \implies \Ext^{i+j}_{D_a^\times}(\mathbf{1},\chi)\]
from Proposition \ref{prop: ext_spectral_seq}. The work has already been done in Corollary \ref{cor: general_ext_sequence} and the previous proposition.

\begin{lemma-N}
Let $r=d^2ef$ and let $\chi$ be a smooth character of $D_a^\times$. There is an exact sequence
\[0 \to \H^1(D_a^\times/I_1, \H^{r-i+1}(I_1,\chi)^\ast) \to \Ext_{D_a^\times}^{i}(\mathbf{1},\chi^\ast) \to (\H^{r-i}(I_1,\chi)^\ast)^{D_a^\times/I_1} \to 0.\]
\end{lemma-N}

\begin{proof}
This follows from Corollary \ref{cor: general_ext_sequence} combined with the previous proposition.
\end{proof}

For $i=r+1$, we get $\Ext^{r+1}_{D_a^\times}(\mathbf{1},\chi^\ast) \simeq \H^1(D_a^\times/I_1, \chi^\ast)$, which is trivial unless $\chi=\mathbf{1}$ in which case we get $\overline{\F}_p$. After $r+1$, the extensions are all trivial.

We can make use of this lemma when $i=r$ as well, since we have results about the surrounding terms in the exact sequence. We have
\begin{display}
0 \ar{r} & \H^1(D_a^\times/I_1, \H^{1}(I_1,\chi)^\ast) \ar{r} & \Ext_{D_a^\times}^{r}(\mathbf{1},\chi^\ast) \ar{r} & (\chi^\ast)^{D_a^\times/I_1} \ar{r} & 0.
\end{display}
We may determine the first group via Proposition \ref{prop:group_cohom_zero} - the dimension is the multiplicity of the trivial representation. This is precisely what was determined in Theorem \ref{thm:char_exts}.

In the case of $d=2$ and $F=\Q_p$, this is enough to determine all extension groups because $r=d^2ef=4$. For a smooth character $\chi$ of $D^\times$, we can compute all $\H^i(I_1,\chi)$ to do this. It is worth noting that this cohomology computation is known and of independent interest in stable homotopy theory, for example see \cite{henn2007finite} Proposition 7. The group $I_1$ appears as the Morava stabilizer group attached to a formal group law, and $\H^\bullet(I_1, \F_p)$ controls certain localization functors $L_{K(n)}$ in Morava $K$-theory. Here, we provide an alternative proof using our results.

We have $\H^0(I_1,\chi)\simeq \chi$, so $\H^4(I_1,\chi) \simeq \chi$. As for $\H^1$, it becomes
\[\H^1(I_1,\chi) \simeq (\overline{\F}_p \oplus \Ind_{D_2^\times}^{D^\times} \chi_{\eta})\otimes \chi,\]
where $\chi_\eta$ is extended from the character of $\F_{p^2}^\times$ sending $x \mapsto x^{p-1}$. The representation $\H^3(I_1,\chi)$ is the dual of $\H^1(I_1,\overline{\F}_p)$ tensored with $\chi$, so it is given by
\[\H^3(I_1,\chi)\simeq (\overline{\F}_p\oplus \Ind_{D_2^\times}^{D^\times} \chi_{\eta}^\ast)\otimes\chi. \]

What remains is to compute $\H^2(I_1,\overline{\F}_p)$. In this case, the Euler characteristic of $I_1$ will be $0$ using the main result of \cite{totaro1999euler}, so we can deduce that $\dim \H^2(I_1,\overline{\F}_p)=4$. Via Lemma \ref{lem:subgrp_generators} and the K\"{u}nneth formula for profinite groups,
\[\H^2(I_1,\overline{\F}_p) \simeq \bigoplus_{i+j=2} \H^i(1+p\Z_p,\overline{\F}_p)\otimes \H^j(I_{1,\Nrd=1},\overline{\F}_p).\]
The isomorphism is via the cross product, which makes this an isomorphism on the level of representations as well because $1+p\Z_p$ and $I_{1,\Nrd=1}$ are normal. As $1+p\Z_p \simeq \Z_p$, its second cohomology is trivial, so there are only two nontrivial terms in the K\"{u}nneth formula. The first is $\H^1(1+p\Z_p,\overline{\F}_p)\otimes \H^1(I_{1,\Nrd=1},\overline{\F}_p)$, which we already know is $\Ind_{D_2^\times}^{D^\times} \chi_{\eta}$ via the decomposition of $\H^1(I_1,\overline{\F}_p)$ in Theorem \ref{thm:homomorphism_classification}. The other nontrivial component is $\H^0(1+p\Z_p,\overline{\F}_p)\otimes \H^2(I_{1,\Nrd=1},\overline{\F}_p)$, which can be found as the image of the Bockstein 
\[\beta_p: \H^1(I_{1,\Nrd=1},\F_p) \to \H^2(I_{1,\Nrd=1},\F_p),\]
after noting that the elements in $\H^1(I_{1,\Nrd=1},\Z/p^2\Z) \simeq \Hom(k_D, \Z/p^2\Z)$ are sent to $0$ so that the map becomes an injection. Note that the isomorphism 
\[\H^1(I_{1,\Nrd=1},\Z/p^2\Z) \simeq \Hom(k_D, \Z/p^2\Z)\]
still follows from the commutator calculation: $[I_{1,\Nrd=1},I_{1,\Nrd=1}]=I_{2,\Nrd=1}$. Tensoring with $\overline{\F}_p$, we see that 
\[\H^2(I_1,\chi) \simeq (\Ind_{D_2^\times}^{D^\times} \chi_{\eta})^{\oplus 2} \otimes \chi.\]
This has dimension four, so we are done. Via Corollary \ref{cor: general_ext_sequence} we can determine all extension groups as well.

\begin{remark}
When $F$ is a more general $p$-adic field or we consider a general division algebra, this method only provides partial information about $\H^2(I_1,\chi)$. The K\"{u}nneth decomposition can provide a subrepresentation of $\H^2(I_1,\overline{\F}_p)$ by considering the image of the Bockstein in $\H^2(I_{1,\Nrd=1},\overline{\F}_p)$ isomorphic to $\H^1(I_{1,\Nrd=1},\overline{\F}_p)$. The other two components are completely determined, as $\H^1(I_1,\overline{\F}_p)$ is known and all cohomology groups of $1+\pi_F \O_F$ are easy to find since it is isomorphic to a direct product of $ef$ copies of $\Z_p$ as a topological group. We can thus use another K\"{u}nneth decomposition on it to show that $\H^i(1+\pi_F \O_F, \overline{\F}_p) \simeq \overline{\F}_p^{\binom{ef}{i}}$. In summary, we have
\[\H^1(I_{1,\Nrd=1},\overline{\F}_p)^{\oplus 2} \oplus \overline{\F}_p^{\binom{ef}{2}} \subset \H^2(I_{1},\overline{\F}_p),\]
where each summand is also as representations of $D^\times/I_1$, with $\overline{\F}_p^{\binom{ef}{2}}$ having the trivial action.
\end{remark}

\section{Acknowledgements}

The authors would like to thank Tasho Kaletha and Karol Koziol for providing valuable mentorship throughout this project, helping simplify several components of this paper, and pointing out errors. They would also like to thank the University of Michigan Math REU for making this project possible. This research was supported by NSF grant DMS-1801687.


\begin{thebibliography}{Wan50}

\bibitem[Bre10]{breuil2010emerging}
Christophe Breuil, \emph{The emerging {$p$}-adic {L}anglands programme},
  Proceedings of the {I}nternational {C}ongress of {M}athematicians. {V}olume
  {II}, Hindustan Book Agency, New Delhi, 2010, pp.~203--230. \MR{2827792}

\bibitem[Hen07]{henn2007finite}
Hans-Werner Henn, \emph{On finite resolutions of {$K(n)$}-local spheres},
  Elliptic cohomology, London Math. Soc. Lecture Note Ser., vol. 342, Cambridge
  Univ. Press, Cambridge, 2007, pp.~122--169. \MR{2330511}

\bibitem[Koz18]{koziol2018hecke}
Karol Kozio\l, \emph{Hecke module structure on first and top
  pro-{$p$}-{I}wahori cohomology}, Acta Arith. \textbf{186} (2018), no.~4,
  349--376. \MR{3879398}

\bibitem[KS21]{koziol2020functorial}
Karol Koziol and David Schwein, \emph{On mod $p$ orientation characters},
  \url{http://www-personal.umich.edu/~kkoziol/orientation.pdf}, 2021.

\bibitem[Laz65]{lazard1965groupes}
Michel Lazard, \emph{Groupes analytiques {$p$}-adiques}, Inst. Hautes
  \'{E}tudes Sci. Publ. Math. (1965), no.~26, 389--603. \MR{209286}

\bibitem[Ly13]{ly2013representations}
Tony Ly, \emph{Repr{\'e}sentations modulo p de {GL}(m, {D}), {D} alg{\`e}bre
  {\`a} division sur un corps local}, Ph.D. thesis, Ph. D. thesis,
  Universit{\'e} Paris Diderot 7, 2013.

\bibitem[Pa{\v{s}}10]{paskunas2007extensions}
Vytautas Pa{\v{s}}k\={u}nas, \emph{Extensions for supersingular representations
  of {${\GL}_2(\Q_p)$}}, Ast\'{e}risque (2010), no.~331, 317--353. \MR{2667891}

\bibitem[PR94]{platonov1993algebraic}
Vladimir Platonov and Andrei Rapinchuk, \emph{Algebraic groups and number
  theory}, Pure and Applied Mathematics, vol. 139, Academic Press, Inc.,
  Boston, MA, 1994, Translated from the 1991 Russian original by Rachel Rowen.
  \MR{1278263}

\bibitem[Rie70]{riehm1970norm}
Carl Riehm, \emph{The norm 1 group of p-adic division algebra}, American
  Journal of mathematics \textbf{92} (1970), no.~2, 499--523.

\bibitem[Sch18]{scholze2015p}
Peter Scholze, \emph{On the {$p$}-adic cohomology of the {L}ubin-{T}ate tower},
  Ann. Sci. \'{E}c. Norm. Sup\'{e}r. (4) \textbf{51} (2018), no.~4, 811--863,
  With an appendix by Michael Rapoport. \MR{3861564}

\bibitem[Ser77]{serre1977linear}
Jean-Pierre Serre, \emph{Linear representations of finite groups}, Graduate
  Texts in Mathematics, Vol. 42, Springer-Verlag, New York-Heidelberg, 1977,
  Translated from the second French edition by Leonard L. Scott. \MR{0450380}

\bibitem[SK79]{shioda1979fermat}
Tetsuji Shioda and Toshiyuki Katsura, \emph{On {F}ermat varieties}, Tohoku
  Math. J. (2) \textbf{31} (1979), no.~1, 97--115. \MR{526513}

\bibitem[Tot99]{totaro1999euler}
Burt Totaro, \emph{Euler characteristics for {$p$}-adic {L}ie groups}, Inst.
  Hautes \'{E}tudes Sci. Publ. Math. (1999), no.~90, 169--225 (2001).
  \MR{1813226}

\bibitem[Wan50]{wang1950commutator}
Shianghaw Wang, \emph{On the commutator group of a simple algebra}, Amer. J.
  Math. \textbf{72} (1950), 323--334. \MR{34380}

\bibitem[Wei49]{weil1949numbers}
Andr\'{e} Weil, \emph{Numbers of solutions of equations in finite fields},
  Bull. Amer. Math. Soc. \textbf{55} (1949), 497--508. \MR{29393}

\end{thebibliography}

\providecommand{\bysame}{\leavevmode\hbox to3em{\hrulefill}\thinspace}
\providecommand{\MR}{\relax\ifhmode\unskip\space\fi MR }
\providecommand{\MRhref}[2]{%
  \href{http://www.ams.org/mathscinet-getitem?mr=#1}{#2}
}
\providecommand{\href}[2]{#2}

\end{document}